\documentclass[10pt,journal,twoside,web]{ieeecolor}
\usepackage{generic}
\usepackage{cite}
\usepackage{amsmath,amssymb,amsfonts}
\usepackage{mathtools}
\mathtoolsset{showonlyrefs}
\usepackage{algorithmic}
\usepackage{graphicx}
\usepackage{textcomp}
\usepackage{subcaption}
\usepackage{makecell}
\usepackage{url}
\newtheorem{remark}{\bf Remark}
\newtheorem{theorem}{\bf Theorem}
\newtheorem{lemma}{\bf Lemma}

\newtheorem{proposition}{\bf Proposition}

\newtheorem{definition}{\bf Definition}

\def\BibTeX{{\rm B\kern-.05em{\sc i\kern-.025em b}\kern-.08em
    T\kern-.1667em\lower.7ex\hbox{E}\kern-.125emX}}
\begin{document}
\title{Distributed Time-Varying Optimization via Unbiased Extremum Seeking}
\author{Xuebin Li, Xuefei Yang, Emilia Fridman, \IEEEmembership{Fellow, IEEE}, Mamadou Diagne, \IEEEmembership{Senior Member, IEEE} and Jiebao Sun
\thanks{This work was supported by the National Natural Science Foundation of China (Grant Nos. 62573155, 62125303, 62188101), Israel Science Foundation (grant no. 446/24), and the Chana and Heinrich Manderman Chair on System Control at Tel Aviv University. Corresponding author: Xuefei Yang.}%
\thanks{Xuebin Li and Jiebao Sun are with the School of Mathematics, Harbin Institute of Technology, China (e-mail: 23B912027@stu.hit.edu.cn; sunjiebao@hit.edu.cn).}%
\thanks{Xuefei Yang is with the Center for Control Theory and Guidance Technology, Harbin Institute of Technology, China (e-mail: xfyang1989@163.com).}%
\thanks{Emilia Fridman is with the School of Electrical and Computer Engineering, Tel-Aviv University, Israel (e-mail: emilia@tauex.tau.ac.il).}%
\thanks{Mamadou Diagne is with the Department of Mechanical and Aerospace Engineering, University of California San Diego, CA 92093, USA (e-mail: mdiagne@ucsd.edu).}%
}

\maketitle

\begin{abstract}
This paper proposes novel distributed optimization algorithms that
address time-varying optimization problems without requiring explicit
gradient information of the objective functions. Traditional distributed methods often rely
on gradient computations, limiting their applicability when only
real-time objective function measurements are available. Leveraging
extremum seeking, we develop continuous-time algorithms
that utilize local measurements and neighbor-shared data to collaboratively
track time-varying optima. Key advancements include compatibility
with directed communication graphs, customizable convergence rates
(asymptotic, exponential, or prescribed-time), and the ability to
handle dynamically evolving objectives. By integrating chirpy probing
signals with time-varying frequencies, our unified framework achieves
accelerated convergence while maintaining stability under mild assumptions.
Theoretical guarantees are established through Lie bracket averaging
and Lyapunov-based analysis, with linear matrix inequality conditions
ensuring rigorous convergence. Numerical simulations validate the
effectiveness of the algorithms.
\end{abstract}


\noindent\emph{Note:} This is the extended version of the manuscript submitted
to the IEEE Transactions on Control of Network Systems. It is identical to the
submitted version except for two additions that the page limit ruled out there.
Appendix \ref{sec:Proof-of-Lemmabound} reports the proof of Lemma
\ref{lem:bound} in full, whereas the submitted version states only the two
bounds that the analysis uses and refers to this document for the intermediate
computations; and Table \ref{tab:notation} at the end of Section
\ref{sec:Preliminaries}-A summarizes the notation. Because both additions shift
the counters, the equation and appendix numbers of this document do not coincide
with those of the submitted version.\par\medskip

\section{Introduction}

Distributed optimization is concerned with a network of agents that
work collaboratively to optimize the sum of their individual objective
functions. This type of problem arises in a variety of applications
in science and engineering, such as optimal power dispatch
in smart grids\cite{ye2016distributed}, distributed machine learning\cite{boyd2011distributed}
and formation control problems \cite{de2023predefined}. {Existing distributed
optimization methods can be formulated in discrete- or continuous-time.} Discrete-time approaches encompass
frameworks such as subgradient methods \cite{nedic2009distributed},
alternating direction methods of multipliers \cite{boyd2011distributed,hong2017distributed,pan2024asynchronous},
exact first-order algorithms \cite{shi2015extra}, and deep-learning based approach \cite{brumali2026data}. Continuous-time
counterparts include distributed proportional-integral-based algorithms
\cite{wang2010control,kia2015distributed}, zero-gradient-sum algorithms
\cite{lu2012zero}, distributed subgradient method \cite{zhu2018continuous},
and distributed optimization with communication delays \cite{yang2016distributed}.
These foundational works generally assume static objective functions,
limiting their applicability to scenarios with time-varying cost functions
such as economic dispatch problems \cite{cherukuri2016initialization}
and cooperative consensus of the unmanned aerial vehicle swarm \cite{yang2024distributed},
where objective functions evolve dynamically.

For the dynamic environments, the optimization objective shifts
from merely identifying static optimal points to dynamically tracking
the trajectory of time-varying optima. One significant approach involves
prediction-correction methodologies \cite{simonetto2016class,fazlyab2017prediction},
which introduced temporal gradient compensation mechanisms to extend
classical Newton approaches to dynamic environments. For scenarios with constraints and limited gradient information,
Chen et al. in \cite{chen2025continuous}
addressed constrained time-varying optimization via incorporating projection maps into continuous-time
zeroth-order dynamics using direct cost function measurements. Building upon the prediction-correction framework,
distributed extensions were developed to enable cooperative strategies for
multi-agent systems \cite{rahili2016distributed,sun2022distributed}. This method was subsequently extended to scenarios with disturbances in \cite{zhang2023continuous}. Wang et al. in \cite{wang2020distributed} established
foundational continuous-time algorithms for resource allocation with
time-varying quadratic costs, which are applicable to both homogeneous time-varying
and heterogeneous time-invariant Hessian configurations. Later, this method was extended
to fully heterogeneous time-varying Hessian scenarios in \cite{wang2022distributed}.
Recent advances further extend the applicability of distributed time-varying optimization: for instance,
Jiang et al. in \cite{jiang2024distributed} eliminated the positive definiteness
requirement for local Hessians, while He et al. in \cite{he2025distributed} demonstrated
fixed-time convergence in constrained optimization through zeroing
neural network implementations.

All aforementioned algorithms fundamentally rely on gradient information
of the objective functions, particularly temporal derivatives of the
gradients. However, in many practical implementations, gradient information
is often unavailable because time-varying objective functions typically lack explicit
 formulations, with only real-time objective value measurements accessible. This intrinsic
dependency renders existing gradient-based methodologies inapplicable
to scenarios where only real-time function measurements are accessible.
To address this issue, an effective approach is to introduce extremum
seeking (ES)--a model-free optimization technique--into distributed
frameworks. ES eliminates the need for explicit mathematical models
by utilizing only instantaneous objective function measurements, thereby
avoiding explicit gradient computations \cite{krstic2000stability,durr2013lie,scheinker2024100}.
Building on this principle, Ye et al. in \cite{ye2016distributed} demonstrated
the integration of ES with saddle-point methods over undirected graphs
for distributed optimization, enabling agents to locate extrema through
measurement data alone without explicit objective function formulations.
Subsequently, Dougherty et al. in \cite{dougherty2016extremum}
developed an ES framework incorporating average consensus protocols
to solve constrained optimization problems with strongly convex social
objectives, though not for consensus-based distributed optimization. Li et al. in \cite{li2020cooperative} formulated cooperative source seeking via networked multi-vehicle systems as a distributed optimization problem and employed stochastic ES to address this challenge.
{More recently, Mimmo et al. in \cite{mimmo2024extremum} combined discrete-time ES with the gradient tracking algorithm \cite{nedic2017achieving} for distributed optimization.}
However, these methodologies remain confined to time-invariant
ES with static objectives and practical stability guarantees, {i.e., convergence to a neighborhood of the extremum rather than exact convergence.}

To address the challenges of time-varying distributed optimization
without explicit gradient information of the objective functions, this paper develops continuous-time distributed
optimization algorithms leveraging unbiased extremum seeking (uES) techniques \cite{yilmaz2023exponential,yilmaz2024asymptotic}.
We apply ES principles to estimate gradient information through
real-time measurements. This gradient estimation mechanism integrates
with proportional-integral-based coordination frameworks, enabling
distributed optimization without explicit gradient computations. For
constant-frequency probing design, we apply Lie bracket averaging
techniques to transform the original error dynamics into a simplified
averaged system. By constructing a Lyapunov function candidate and
performing rigorous stability analysis, we derive verifiable linear matrix inequality (LMI) conditions
that guarantee stability of the averaged system given initial values
within an invariant set. This framework eventually establishes uniform
asymptotic stability of the original distributed ES
system. To further enhance convergence rates beyond asymptotic guarantees (exponential, prescribed-time),
we extend this framework by adopting the chirpy probing paradigm. By applying a time-scale transformation
to the averaged system, we convert the analysis of chirpy probing
cases into equivalent constant-frequency probing problems for unified
treatment.

The contributions of this paper are summarized as follows:
\begin{enumerate}
\item Novel distributed optimization algorithm: We propose a class of distributed
optimization algorithms that require neither explicit gradient information
nor prior knowledge of cost functions. Each agent exclusively
utilizes real-time local measurements and neighbor-shared information (or measured relative state information)
for collaborative optimization. The framework can solve time-varying
distributed optimization problems under mild assumptions on the time-varying
objective functions.
\item Advancements over existing distributed ES approaches: Compared with
prior distributed ES methods \cite{ye2016distributed}{,\cite{mimmo2024extremum}},
our framework demonstrates three key advancements:
\begin{enumerate}
\item Directed graph compatibility: Our algorithm can operate over directed
graphs, whereas the continuous-time algorithm in \cite{ye2016distributed} is restricted to undirected
communication topologies.
\item Unbiased convergence with rate customization: We achieve unbiased
convergence with customizable rates (asymptotic, exponential, or prescribed-time),
while \cite{ye2016distributed}{,\cite{mimmo2024extremum}} only guarantee practical exponential-stability.
\item Time-varying function handling: The proposed methodology is applicable to
time-varying optimization problems with dynamic objectives, a capability
absent in \cite{ye2016distributed}{,\cite{mimmo2024extremum}}.
\end{enumerate}
\end{enumerate}
{It is worth noting that extending the uES framework 
from a single-agent setting to the multi-agent distributed 
setting is nontrivial and poses two critical challenges.} 
First, the optimal solution of the aggregate objective 
(the sum of each agent's objective function) does not 
necessarily coincide with the local minimizers of 
individual agents' objective functions. This discrepancy 
leads to unfavorable properties of local cost functions 
near the global optimum (e.g., non-vanishing gradients), 
which significantly complicates both algorithm design and 
stability analysis. Second, the averaged system obtained 
via Lie bracket averaging exhibits substantially greater 
complexity than that in the single-agent case, 
leading to a more challenging stability analysis. {For 
more details, we refer to 
Remarks~\ref{rem:averaged_system} and ~\ref{rem:single_vs_distributed}.}\par

The rest of the paper is organized as follows. Section \ref{sec:Preliminaries}
introduces some basic concepts and useful lemmas. Section \ref{sec:Problem-Definition}
outlines the problem formulation. Section \ref{subsec:constant-f}
introduces a constant-frequency probing scheme with asymptotic convergence
guarantees. Section \ref{subsec:chirpy} extends this framework through
a unified formulation employing chirpy probing signals, which enables
asymptotic, exponential, and prescribed-time convergence via choosing
the corresponding growth functions. Section \ref{subsec:Numerical-Simulation}
presents numerical simulations to validate these algorithms' effectiveness
and comparative performance. Finally, Section \ref{sec:Conclusion}
concludes the paper.

\section{Preliminaries}

\label{sec:Preliminaries}

This section introduces mathematical notation and foundational concepts
from convex analysis and graph theory.

\subsection{Definitions and Concepts}

Let $\mathbb{R}$ denote the real numbers, $\mathbf{1}_{N}$ (resp.
$\mathbf{0}_{N}$) the $N$-dimensional ones (resp. zeros) vector,
and $I_{N}$ the $N\times N$ identity matrix. For $s=1,\dots,d$, let
$e_{s}\in\mathbb{R}^{d}$ denote the $s$-th standard basis vector whose $s$-th
entry is $1$ and the others are $0$. For matrices $A\in\mathbb{R}^{n\times m}$
and $B\in\mathbb{R}^{p\times q}$, we let $A\otimes B$ denote their
Kronecker product. The notation $\|\cdot\|$ refers to both the usual
Euclidean norm for vectors and the induced 2-norm for matrices. Derivatives
follow the following layout convention: for scalar-valued $f:\mathbb{R}^{d}\rightarrow\mathbb{R},\nabla f(x)=\frac{\partial}{\partial x}f(x)$
is a column vector, and the Jacobian of $g:\mathbb{R}^{m}\rightarrow\mathbb{R}^{n}$
is structured as $\frac{\partial g}{\partial x}\in\mathbb{R}^{n\times m}$.
The Lie bracket of two vector fields $f,g:\mathbb{R}^{n}\times\mathbb{R}\rightarrow\mathbb{R}^{n}$
with $f(\cdot,t),g(\cdot,t)$ being continuously differentiable is
defined by $[f,g](x,t):=\frac{\partial g(x,t)}{\partial x}f(x,t)-\frac{\partial f(x,t)}{\partial x}g(x,t)$.
A weighted directed graph $\mathcal{G}=(\mathcal{V},\mathcal{E},\mathcal{A})$
consists of nodes $\mathcal{V}=\{1,\ldots,N\}$, edges $\mathcal{E}\subseteq\mathcal{V}\times\mathcal{V}$,
and an adjacency matrix $\mathcal{A}$. $(i,j)\in\mathcal{E}$ means
that agent $i$ can receive information from agent $j$. The adjacency
matrix $\mathcal{A}$ has entries $a_{ij}=1$ if $(i,j)\in\mathcal{E}$
(0 otherwise, with $a_{ii}\equiv0$ ). A directed path is a sequence
of nodes connected by edges. A digraph is strongly connected if for
every pair of nodes there is a directed path connecting them. The
out-degree and in-degree of a node $i$ are $d_{\text{out }}^{i}=\sum_{j=1}^{N}a_{ij}$
and $d_{\text{in }}^{i}=\sum_{j=1}^{N}a_{ji}$, respectively. A graph
$\mathcal{G}$ is weight-balanced if $d_{\text{in }}^{i}=d_{\text{out }}^{i}$
for any node $i\in\mathcal{V}$. The Laplace matrix of $\mathcal{G}$
is $L=D_{\text{out }}-\mathcal{A}$, where $D_{\text{out }}=\operatorname{diag}\left\{ d_{\text{out }}^{1},d_{\text{out }}^{2},\ldots,d_{\text{out }}^{N}\right\} \in\mathbb{R}^{N\times N}$.
If $\mathcal{G}$ is strongly connected, then $L$ has a simple zero
eigenvalue with eigenvector $\mathbf{1}_{N}$. For multi-agent systems,
we extend the Laplacian as $\mathbf{L}=L\otimes I_{d}$ to deal with
variables of dimension $d\in \mathbb{N}$.

Let $x=\left(x_{1},\dots,x_{d}\right)^{\mathrm{T}}\in\mathbb{R}^{d}$
and $a\in\mathbb{R}$, for notational conciseness, we adopt the following
symbolic conventions: $\cos(x)=\mathrm{col}\{\cos(x_{s})\}_{s=1}^{d}$, $\sin(x)=\mathrm{col}\{\sin(x_{s})\}_{s=1}^{d}$,
$a+x=\mathrm{col}\{a+x_{s}\}_{s=1}^{d}$.

{The main notation used throughout the paper is summarized in Table~\ref{tab:notation}.}

\begin{table}[!t]
\centering
\caption{Summary of Main Notation}
\label{tab:notation}
\small
\setlength{\tabcolsep}{3pt}
\begin{tabular}{@{}l p{0.60\columnwidth}@{}}
\hline
Symbol & Meaning \\
\hline
\multicolumn{2}{@{}l}{\textit{Graphs, vectors, and operators}}\\
$N,\ d,\ l$ & number of agents; dimension of $x$; dimension of $\zeta$ \\
$\mathcal{G},\mathcal{V},\mathcal{E},\mathcal{A}$ & digraph, node set, edge set, adjacency matrix \\
$L,\ \mathbf{L}$ & Laplacian $L=D_{\mathrm{out}}-\mathcal{A}$; $\mathbf{L}=L\otimes I_{d}$ \\
$\mathcal{L}$ & reduced Laplacian $\mathcal{L}=R^{\mathrm{T}}LR$ (appendix) \\
$r,\ R$ & $r=\tfrac{1}{\sqrt{N}}\mathbf{1}_{N}$; $R^{\mathrm{T}}r=\mathbf{0}$, $R^{\mathrm{T}}R=I_{N-1}$ \\
$\mathbf{1}_{N},\mathbf{0}_{N},I_{N},e_{s}$ & ones vector, zeros vector, identity, $i$th basis vector \\
$\otimes,\ \nabla,\ [\cdot,\cdot]$ & Kronecker product, gradient, Lie bracket \\
$\operatorname{col}\{\cdot\},\ \|\cdot\|$ & column stack, Euclidean / induced $2$-norm \\
\hline
\multicolumn{2}{@{}l}{\textit{Signals and variables}}\\
$f_{i},\ f$ & local cost of agent $i$; total cost $\sum_{i}f_{i}$ \\
$\zeta(t)$ & unknown exogenous signal, $\zeta:\mathbb{R}\to\mathbb{R}^{l}$ \\
$x_{i}$ & decision variable (control input) of agent $i$ \\
$x^{*}(t),\ x_{i}^{*}(t)$ & optimizer of $f$; minimizer of local $f_{i}$ \\
$\mathbf{z},\ \eta$ & auxiliary (integral) variable; high-pass filter state \\
$m_{i},\ M_{i}$ & strong-convexity and gradient-Lipschitz constants \\
\hline
\multicolumn{2}{@{}l}{\textit{Extremum seeking and time-scaling parameters}}\\
$\omega,\ \omega_{s}$ & probing frequency; distinct probing frequencies \\
$\alpha,\ k,\ \omega_{h}$ & probing amplitude; gain; high-pass filter cutoff \\
$\xi(t)$ & time-scaling in Sec.~\ref{subsec:constant-f}, $\xi=(1+\beta(t-t_{0}))^{1/v}$ \\
$\phi(t)$ & time-scaling (chirpy probing) in Sec.~\ref{subsec:chirpy} \\
$\beta,v,p,q,\rho$ & parameters of the time-scaling functions \\
$c$ & decay-rate exponent in Assumption~2 \\
\hline
\end{tabular}
\end{table}

\subsection{Useful Lemmas}
\begin{lemma}
\cite{boyd2004convex} \label{lem:Lipschitz_grad}The following properties
are equivalent for a twice-differentiable function $f:\mathbb{R}^{d}\rightarrow\mathbb{R}$:

1. $\nabla f(x)$ is M-Lipschitz on $\mathrm{\mathbb{R}}^{d}$.

2. $(\nabla f(x)-\nabla f(y))^{\mathrm{T}}(x-y)\leq M\|x-y\|^{2},\forall x,y\in\mathbb{R}^{d}$.

3. $\nabla^{2}f(x) \leq MI,\forall x\in\mathbb{R}^{d}$.

4. $f(y)\leq f(x)+\nabla f(x)^{\mathrm{T}}(y-x)+\frac{M}{2}\|y-x\|^{2}$.
\end{lemma}
\begin{lemma}
\cite{boyd2004convex} \label{lem:strong_convex}The following properties
are equivalent for a twice-differentiable function $f:\mathbb{R}^{d}\rightarrow\mathbb{R}$:

1. $f$ is $m$-strongly convex on \textup{$\mathrm{\mathbb{R}}^{d}$}.

2. $(\nabla f(x)-\nabla f(y))^{\mathrm{T}}(x-y)\geq m\|x-y\|^{2},\forall x,y\in\mathbb{R}^{d}$.

3. $\nabla^{2}f(x) \geq mI,\forall x\in\mathbb{R}^{d}$.

4. $f(y)\geq f(x)+\nabla f(x)^{\mathrm{T}}(y-x)+\frac{m}{2}\|y-x\|^{2},\forall x,y\in\mathbb{R}^{d}$.
\end{lemma}

\section{Problem Definition}

\label{sec:Problem-Definition}

Consider a network of $N$ agents indexed by $\mathcal{V}=\{1,\ldots,N\}$,
communicating through a graph $\mathcal{G}=(\mathcal{V},\mathcal{E},\mathcal{A}).$
 We aim to minimize the optimization objective
\begin{equation*}
 f(x,\zeta(t))=\sum_{i=1}^{N}f_{i}(x,\zeta(t)),
\end{equation*}
where $x\in\mathbb{R}^{d}$ represents the decision variable, $\zeta:\mathbb{R}\rightarrow\mathbb{R}^{l}$
is an unknown function, $f_{i}:\mathrm{\mathbb{R}}^{d}\times\mathbb{R}^{l}\rightarrow\mathrm{\mathbb{R}}$
represents the local cost function for agent $i$. The explicit expressions
of $f_{i}(x,\zeta(t))$, and their gradients are unknown. Only
measurements of the function $f_{i}(x,\zeta(t))$ are available
to agent $i$.

We make several assumptions. First, we establish the structural properties of the cost functions:

\textbf{Assumption 1} For each agent $i\in\mathcal{V}$
and time $t\in\mathbb{R}$, the following holds, where the constants
$m_{i}$ and $M_{i}$ are independent of $t$ and hence of $\zeta(t)$:

1. $f_{i}(\cdot,\cdot)$ is three times continuously differentiable.

2. $f_{i}(\cdot,\zeta(t))$ is $m_{i}$-strongly convex.

3. $\frac{\partial}{\partial x}f_{i}(\cdot,\zeta(t))$ is $M_{i}$-Lipschitz
continuous.

Under these conditions, the objective function $f(x,\zeta(t))$
inherits strong convexity, which guarantees that there exists a unique
continuously differentiable solution mapping $\pi:\mathbb{R}^{l}\rightarrow\mathbb{R}^{d}$
such that:
\[
x^{*}(t)=\pi\left(\zeta(t)\right),
\]
where $x^{*}(t)$ denotes the time-varying optimal solution. Each
local function $f_{i}$ is strongly convex with its own minimizer
$x_i^{*}(t)\in\mathbb{R}^{d}$. Similar assumptions can be found in
\cite{gharesifard2013distributed}, \cite{kia2015distributed} and \cite{yang2016distributed}.

Next, we assume bounds on $\zeta(t)$, $x^{*}(t)$, $x_i^{*}(t)$ and their derivatives:

\textbf{Assumption 2.} The time-varying functions
satisfy:
\[
\begin{aligned}
\left\Vert \zeta(t)\right\Vert +\left\Vert x^{*}(t)\right\Vert +\sum_{i=1}^{N}\left\Vert x_i^{*}(t)\right\Vert  &\leq M_{1},\\
\left\Vert \dot{\zeta}(t)\right\Vert +\left\Vert \ddot{\zeta}(t)\right\Vert +\left\Vert \dot{x}^{*}(t)\right\Vert +\left\Vert \ddot{x}^{*}(t)\right\Vert  &\leq M_{2}\phi^{c}(t),
\end{aligned}
\]
where $M_{1},M_{2}\geq0$ are unknown constants, the power $c\in\mathbb{R}$,
and $\phi:\mathbb{R}\rightarrow\mathbb{R}$ is a strictly increasing
function.

Assumption 2 restricts both the range of the time-varying functions
and their rate of change. For detailed discussions on the parameter $c$, we refer the reader to Remark \ref{remark:c<0} and Remark \ref{remark:cp}.

Finally, we specify the network connectivity requirements:

\textbf{Assumption 3.} The communication graph
$\mathcal{G}$ is weight-balanced and strongly connected.

Under Assumptions 1--3, our objective is to design the inputs $x_{i}(t)$ such that they collectively drive the system to the optimal solution $x^{*}(t)$ that minimizes $f(x,\zeta(t))$.

\section{Distributed Optimization Algorithms via uES with Constant-Frequency Probing}

\label{subsec:constant-f}

\begin{figure}[!t]
\centering
\includegraphics[width=\columnwidth]{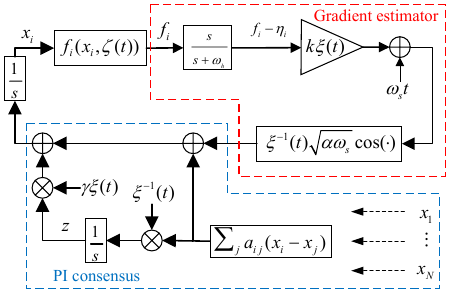}
\caption{Block diagram of the proposed distributed unbiased extremum-seeking architecture, shown for agent $i$.}
\label{fig:blockdiagram}
\end{figure}

Let $\mathbf{x}=\left(x_{1}^{\mathrm{T}},\ldots,x_{N}^{\mathrm{T}}\right)^{\mathrm{T}}\in\left(\mathbb{R}^{d}\right)^{N}$
represent the concatenated state vector where $x_{i}$ denotes agent
$i$'s estimate of the optimal solution $x^{*}$. To streamline subsequent
analysis, we introduce the following formal notations: 
\begin{equation*}
\eta = \mathrm{col}\{\eta_{i}\}_{i=1}^{N}, \quad h(\mathbf{x},\zeta(t)) = \mathrm{col}\{f_{i}(x_{i},\zeta(t))\}_{i=1}^{N},
\end{equation*}
from which it follows that
\begin{equation*}
h\left(\boldsymbol{1}_{N}\otimes x^{*}(t),\zeta(t)\right)=\mathrm{col}\{f_{i}(x^{*}(t),\zeta(t))\}_{i=1}^{N}.
\end{equation*}
Inspired by \cite{kia2015distributed} and \cite{yilmaz2024asymptotic},
we propose the following distributed optimization algorithm (shown in Figure \ref{fig:blockdiagram}).
For each agent $i\in\mathcal{V}$, the local update law is
\begin{equation*}
\begin{aligned}
\dot{x}_{i} &= \xi^{-1}(t)\sum_{s=1}^{d}\sqrt{\alpha\omega_{s}}
              \cos\!\big(\omega_{s}t+k\xi(t)(f_{i}(x_{i},\zeta(t))-\eta_{i})\big)e_{s}\\
            &\quad -\sum_{j=1}^{N}a_{ij}(x_{i}-x_{j})-z_{i}\xi^{-1}(t),\\
\dot{\eta}_{i} &= -\omega_{h}\eta_{i}+\omega_{h}f_{i}(x_{i},\zeta(t)),\\
\dot{z}_{i} &= \gamma\sum_{j=1}^{N}a_{ij}(x_{i}-x_{j})\,\xi(t).
\end{aligned}
\end{equation*}
where $a_{ij}$ are the entries of the adjacency matrix $\mathcal{A}$, \{$\alpha,\beta,\gamma,k,v$\} are positive parameters, the frequencies $\omega_{s}=\omega\hat{\omega}_{s}$ are chosen such
that $\hat{\omega}_{s}\in\mathbb{N}$ and $\hat{\omega}_{s}\neq\hat{\omega}_{\ell},\forall s\neq \ell$, and the asymptotically growing function $\xi(t)$ is given by
\begin{equation}
\xi(t)=\left(1+\beta\left(t-t_{0}\right)\right)^{\frac{1}{v}}.\label{eq:xi}
\end{equation}
Equivalently, in the compact network-level notation $\mathbf{x}=\operatorname{col}\{x_{i}\}_{i=1}^{N}$, $\eta=\operatorname{col}\{\eta_{i}\}_{i=1}^{N}$, $\mathbf{z}=\operatorname{col}\{z_{i}\}_{i=1}^{N}$, the algorithm reads
\begin{equation}
\begin{aligned}
\dot{\mathbf{x}} &= \sum_{s=1}^{d}\xi^{-1}(t)\sqrt{\alpha\omega_{s}}\cos\left(\omega_{s}t+k\xi(t)\left(h(\mathbf{x},\zeta(t))-\eta\right)\right)\\
&\quad \otimes e_{s}+\left(-\mathbf{L}\mathbf{x}-\mathbf{z}\cdot\xi^{-1}(t)\right),\\
\dot{\eta} &= -\omega_{h}\eta+\omega_{h}h(\mathbf{x},\zeta(t)),\\
\dot{\mathbf{z}} &= \gamma\mathbf{L}\mathbf{x}\cdot\xi(t),
\end{aligned}
\label{algo_extremum}
\end{equation}
Here $\cos(\cdot)$ and the scalar-plus-vector addition are componentwise (Section \ref{sec:Preliminaries}), so $\cos(\cdot)\in\mathbb{R}^{N}$ and $\cos(\cdot)\otimes e_{s}\in\mathbb{R}^{Nd}$.
\begin{remark}
This algorithm is distributed because each agent only uses its private
cost function $f_{i}(x_{i},\zeta(t))$ and the information received
from its neighbors regarding their corresponding variables in $\mathbf{x}$.
Each agent utilizes real-time measurements to estimate the gradient
terms (see the first term in the first equation of \eqref{algo_extremum} and the red dashed box in Figure \ref{fig:blockdiagram}, adapted from bounded ES \cite{scheinker2014extremum}), thereby eliminating the dependence on explicit
gradient computations required in model-based approaches.
 The distributed coordination terms---the Laplacian feedback $\mathbf{L}\mathbf{x}$ (proportional) and the auxiliary variable $\mathbf{z}$ (integral)---form a PI consensus mechanism \cite{kia2015distributed,yang2016distributed} (the blue dashed box in Figure \ref{fig:blockdiagram}) that enforces consensus among agents.
The asymptotically growing function $\xi(t)$ is employed
to guarantee unbiased convergence, and the
high-pass filters provide better transient
behavior by removing possible constant offsets of the individual cost
function $f_{i}$ (see \cite{durr2013lie}).
\end{remark}

Next, we analyze the stability of algorithm \eqref{algo_extremum}.
Consider the following transformations
\begin{equation}
\begin{aligned}
\tilde{\mathbf{x}}_{f} &= \xi(t)(\mathbf{x}-\boldsymbol{1}_{N}\otimes x^{*}(t)),\\
\tilde{\eta}_{f} &= \xi(t)(\eta-h(\boldsymbol{1}_{N}\otimes x^{*}(t),\zeta(t))),\\
\tilde{\mathbf{z}}_{f} &= \mathbf{z}-\mathbf{z}^{*}(t),
\end{aligned}
\label{trans}
\end{equation}
where 
\begin{equation}
\mathbf{z}^{*}(t)=-\frac{1}{2}\alpha k\frac{\partial G}{\partial\mathbf{x}}(\boldsymbol{1}_{N}\otimes x^{*}(t),\zeta(t))\label{eq:z_star}
\end{equation}
with
\begin{equation}\label{eq:G}
\frac{\partial G}{\partial\mathbf{x}}(\mathbf{x},\zeta(t))=\operatorname{col}\left\{\frac{\partial f_{i}\left(x_{i},\zeta(t)\right)}{\partial x_{i}}\right\}_{i=1}^{N}.
\end{equation}
{Using the expression of $\xi(t)$ and differentiating \eqref{trans} with respect to time, we convert \eqref{algo_extremum} into the following system:}
\allowdisplaybreaks
\begin{equation}
\begin{aligned}
\dot{\tilde{\mathbf{x}}}_{f}={} & \frac{\beta}{v}\xi^{-v}(t)\tilde{\mathbf{x}}_{f}-\xi(t)\boldsymbol{1}_{N}\otimes\dot{x}^{*}(t)\\
 & \hspace{-1.5em}+\sum_{s=1}^{d}\sqrt{\alpha\omega_{s}}\cos(\omega_{s}t)\cos\left(k\xi(t)h_{f}(\tilde{\mathbf{x}}_{f},\zeta(t))-k\tilde{\eta}_{f}\right)\otimes e_{s}\\
 & \hspace{-1.5em}-\sum_{s=1}^{d}\sqrt{\alpha\omega_{s}}\sin(\omega_{s}t)\sin\left(k\xi(t)h_{f}(\tilde{\mathbf{x}}_{f},\zeta(t))-k\tilde{\eta}_{f}\right)\otimes e_{s}\\
+{} & \left(-\mathbf{L}\tilde{\mathbf{x}}_{f}-\tilde{\mathbf{z}}_{f}-\mathbf{z}^{*}(t)\right),\\
\dot{\tilde{\eta}}_{f}={} & \left(\frac{\beta}{v}\xi^{-v}(t)-\omega_{h}\right)\tilde{\eta}_{f}\\
 & +\omega_{h}\xi(t)h_{f}(\tilde{\mathbf{x}}_{f},\zeta(t))-\xi(t)\dot{h}(\boldsymbol{1}_{N}\otimes x^{*}(t)),\\
\dot{\tilde{\mathbf{z}}}_{f}={} & \gamma\mathbf{L}\mathbf{x}\cdot\xi(t)-\dot{\mathbf{z}}^{*}(t)=\gamma\mathbf{L}\tilde{\mathbf{x}}_{f}-\dot{\mathbf{z}}^{*}(t),
\end{aligned}
\label{eq:error_system}
\end{equation}
where
\begin{equation}
\begin{aligned}
&h_{f}(\tilde{\mathbf{x}}_{f},\zeta(t))=h(\mathbf{x},\zeta(t))-h(\boldsymbol{1}_{N}\otimes x^{*}(t),\zeta(t))\\
&= h(\xi^{-1}(t)\tilde{\mathbf{x}}_{f}+\boldsymbol{1}_{N}\otimes x^{*}(t),\zeta(t))-h(\boldsymbol{1}_{N}\otimes x^{*}(t),\zeta(t)).
\end{aligned}
\label{eq:h_f}
\end{equation}
In the third equation of \eqref{eq:error_system} we used $\mathbf{L}\left(\boldsymbol{1}_{N}\otimes x^{*}(t)\right)=\left(L\boldsymbol{1}_{N}\right)\otimes x^{*}(t)=0$, which gives $\mathbf{L}\mathbf{x}\cdot\xi(t)=\mathbf{L}\tilde{\mathbf{x}}_{f}$.
With these transformations, we only need to demonstrate
the practical stability of the origin of system \eqref{eq:error_system}.
Given the asymptotic growth property of $\xi(t)$, this would ultimately
ensure the asymptotic convergence of $\mathbf{x}$ to $\boldsymbol{1}_{N}\otimes x^{*}$.
To analyze system \eqref{eq:error_system}, we employ the Lie bracket
averaging technique \cite{durr2013lie,durr2013saddle} commonly used in ES. 
In brief, it allows a complex system driven by high-frequency oscillations to be approximated by a relatively simple \emph{averaged system}, which here amounts to a gradient-descent term.
To this end,
we first reformulate system \eqref{eq:error_system} into the input-affine
form. Specifically, we introduce the vector $\varsigma=\left(\tilde{\mathbf{x}}_{f}^{\mathrm{T}},\tilde{\eta}_{f}^{\mathrm{T}},\tilde{\mathbf{z}}_{f}^{\mathrm{T}}\right)^{\mathrm{T}}$,
then system \eqref{eq:error_system} can be rewritten as
\vspace{-0.5em}
\begin{equation}
\begin{aligned}
\dot{\varsigma} &= b_{0}(\varsigma,t)+\sum_{s=1}^{d}b_{1s}(\varsigma,t)\sqrt{\omega_{s}}\cos(\omega_{s}t)\\
&\quad +\sum_{s=1}^{d}b_{2s}(\varsigma,t)\sqrt{\omega_{s}}\sin(\omega_{s}t)
\end{aligned}
\label{rewrite_error}
\end{equation}
\vspace{-0.5em}
with
\[
b_{0}(\varsigma,t)=\begin{bmatrix}\dfrac{\beta}{v}\xi^{-v}(t)\tilde{\mathbf{x}}_{f}-\mathbf{L}\tilde{\mathbf{x}}_{f}-\tilde{\mathbf{z}}_{f}-\mathbf{z}^{*}(t)\\
\quad-\xi(t)\boldsymbol{1}_{N}\otimes\dot{x}^{*}(t)\\[2ex]
\left(\dfrac{\beta}{v}\xi^{-v}(t)-\omega_{h}\right)\tilde{\eta}_{f}+\omega_{h}\xi(t)h_{f}(\tilde{\mathbf{x}}_{f},\zeta(t))\\
\quad-\xi(t)\dot{h}(\boldsymbol{1}_{N}\otimes x^{*}(t),\zeta(t))\\[2ex]
\gamma\mathbf{L}\tilde{\mathbf{x}}_{f}-\dot{\mathbf{z}}^{*}(t)
\end{bmatrix},
\]
\[
b_{1s}(\varsigma,t)=\left[\begin{array}{c}
\sqrt{\alpha}\cos\left(k\xi(t)h_{f}(\tilde{\mathbf{x}}_{f},\zeta(t))-k\tilde{\eta}_{f}\right)\otimes e_{s}\\
\mathbf{0}_{N+Nd}
\end{array}\right],
\]
\[
b_{2s}(\varsigma,t)=\left[\begin{array}{c}
-\sqrt{\alpha}\text{\ensuremath{\sin}}\left(k\xi(t)h_{f}(\tilde{\mathbf{x}}_{f},\zeta(t))-k\tilde{\eta}_{f}\right)\otimes e_{s}\\
\mathbf{0}_{N+Nd}
\end{array}\right].
\]
{Next, we need to verify that the system \eqref{rewrite_error} satisfies the
boundedness assumptions in Assumption A2 of \cite{durr2013lie}, which we
formalize in the following lemma.}
\begin{lemma}
\label{lem:bound}Suppose Assumptions 1--2 hold with $\phi(t)=\xi(t)$ and $c \leq -1$.
Let \textbf{$\mathcal{X}\subset\mathbb{R}^{2Nd+N}$} be a compact
set. The norms
$\|b_{0}(\varsigma,t)\|$, $\|b_{1s}(\varsigma,t)\|$,
$\|b_{2s}(\varsigma,t)\|$, $\|\frac{\partial b_{0}(\varsigma,t)}{\partial\varsigma}\|$, $\|\frac{\partial b_{1s}(\varsigma,t)}{\partial\varsigma}\|$,$\|\frac{\partial b_{2s}(\varsigma,t)}{\partial\varsigma}\|$, $\|\frac{\partial b_{0}(\varsigma,t)}{\partial t}\|$,
$\|\frac{\partial b_{1s}(\varsigma,t)}{\partial t}\|$, $\|\frac{\partial b_{2s}(\varsigma,t)}{\partial t}\|$,\\$\|\frac{\partial\left[b_{1s}(\varsigma,t),b_{2s}(\varsigma,t)\right]}{\partial\varsigma}\|$ and $\|\frac{\partial\left[b_{1s}(\varsigma,t),b_{2s}(\varsigma,t)\right]}{\partial t}\|$
are all bounded for all $\varsigma\in\mathcal{X}$, $t\in\mathbb{R}$,
$s=1,\dots,d$.
\end{lemma}
\begin{proof}
See Appendix \ref{sec:Proof-of-Lemmabound}.
\end{proof}

With the boundedness conditions established in Lemma \ref{lem:bound},
we now perform the Lie bracket averaging operation on system \eqref{rewrite_error}.
Denoting $\bar{\varsigma}=\left(\bar{\mathbf{x}}_{f}^{\mathrm{T}},\bar{\eta}_{f}^{\mathrm{T}},\bar{\mathbf{z}}_{f}^{\mathrm{T}}\right)^{\mathrm{T}}$,
we can derive the Lie bracket averaged system from \eqref{rewrite_error}:
\vspace{-0.5em}
\[
\dot{\bar{\varsigma}}=b_{0}(\bar{\varsigma},t)+\frac{1}{2}\sum_{s=1}^{d}\left[b_{1s}(\bar{\varsigma},t),b_{2s}(\bar{\varsigma},t)\right].
\]
\vspace{-0.5em}
We express it as:
\begin{equation}
\begin{aligned}
\dot{\bar{\mathbf{x}}}_{f}={} & \frac{\beta}{v}\xi^{-v}(t)\bar{\mathbf{x}}_{f}-\frac{1}{2}\alpha k\frac{\partial G}{\partial\mathbf{x}}(\bar{\mathbf{x}},\zeta(t))\\
 & +\left(-\mathbf{L}\bar{\mathbf{x}}_{f}-\bar{\mathbf{z}}_{f}-\mathbf{z}^{*}(t)\right)-\xi(t)\boldsymbol{1}_{N}\otimes\dot{x}^{*}(t),\\
\dot{\bar{\eta}}_{f}={} & \left(\frac{\beta}{v}\xi^{-v}(t)-\omega_{h}\right)\bar{\eta}_{f}\\
 & +\omega_{h}\xi(t)h_{f}(\bar{\mathbf{x}}_{f},\zeta(t))-\xi(t)\dot{h}(\boldsymbol{1}_{N}\otimes x^{*}(t),\zeta(t)),\\
\dot{\bar{\mathbf{z}}}_{f}={} & \gamma\mathbf{L}\bar{\mathbf{x}}_{f}-\dot{\mathbf{z}}^{*}(t),
\end{aligned}
\label{averaged_system}
\end{equation}
where we have used the notation 
$\bar{\mathbf{x}}=\left[\bar{x}_{1}^{\mathrm{T}},\cdots,\bar{x}_{N}^{\mathrm{T}}\right]^{\mathrm{T}}=\xi^{-1}(t)\bar{\mathbf{x}}_{f}+\boldsymbol{1}_{N}\otimes x^{*}(t).$
{
\begin{remark}\label{rem:averaged_system}
Compared with the model-based algorithm in \cite{kia2015distributed}, the averaged system \eqref{averaged_system} exhibits three key distinctions.
First, the term $\frac{\beta}{v}\xi^{-v}(t)\bar{\mathbf{x}}_{f}$ constitutes a \emph{positive feedback} arising from the scaling transformation \eqref{trans}, introducing a destabilizing effect absent in \cite{kia2015distributed}.
Second, $\frac{\partial G}{\partial\mathbf{x}}(\bar{\mathbf{x}},\zeta(t))$ (defined in \eqref{eq:G}) is not a standard gradient of the error variables. Indeed,
\[
\frac{\partial G}{\partial\mathbf{x}}(\bar{\mathbf{x}},\zeta(t))
= \xi(t)\operatorname{col}_{i=1}^{N}\!\left\{\nabla_{\bar{x}_{fi}}f_{i}\!\left(\xi^{-1}(t)\bar{x}_{fi}+x^{*}(t),\zeta(t)\right)\right\},
\]
which differs structurally from $\nabla f_{i}$ evaluated at the error variables due to the scaling factor $\xi(t)$ and the nonlinear change of variables.
Third, the terms $-\mathbf{z}^{*}(t)$ (defined in \eqref{eq:z_star}), $-\xi(t)\boldsymbol{1}_{N}\otimes\dot{x}^{*}(t)$, and $-\dot{\mathbf{z}}^{*}(t)$ encode the time derivatives of the optimal trajectory and can be regarded as perturbations that vanish as $t\to\infty$ under Assumption~2 (see \eqref{varying}).
On the other hand, compared to the single-agent case, our system~\eqref{averaged_system} is considerably more complex, and consequently the stability analysis is also significantly more involved.
\end{remark}}

{As highlighted in Remark~\ref{rem:averaged_system}, the positive feedback term $\frac{\beta}{v}\xi^{-v}(t)\bar{\mathbf{x}}_{f}$ and the perturbations induced by the time-varying optimal trajectory introduce destabilizing effects into system~\eqref{averaged_system}. To ensure that the stabilizing terms dominate these adverse contributions, the parameter $v$ must be sufficiently large so that the gain $\frac{\beta}{v}\xi^{-v}(t)$ decays rapidly, and the parameter $c$ in Assumption~2 must be sufficiently negative to limit the rate of temporal variation. Specifically, the subsequent analysis requires
\begin{equation}\label{eq:vc_condition}
v\geq2,\qquad c<-3.
\end{equation}}

We now investigate
the stability properties of system \eqref{averaged_system}. Prior to this analysis, we first present
a critical property of system \eqref{averaged_system}. The proof is provided in Appendix
\ref{app_em:equilibrium}:
\begin{lemma}
\label{lem:invariant_set}Let Assumption 1 be satisfied. Given any
$\epsilon\in\mathbb{R}^{d}$, $\mathcal{\varOmega}\left(\epsilon\right)=\left\{ (\mathbf{x},\eta,\mathbf{z})|\left(\mathbf{1}_{N}\otimes I_{d}\right)^{\mathrm{T}}\mathbf{z}=\epsilon\right\} $
is a positive invariant set of system \eqref{averaged_system}.
\end{lemma}
It is the same condition on $\bar{\mathbf{z}}_{f}$, since \eqref{eq:z_star} and optimality of $x^{*}$ give $(\mathbf{1}_{N}\otimes I_{d})^{\mathrm{T}}\mathbf{z}^{*}=\mathbf{0}_{d}$.
The stability conclusions for the averaged system \eqref{averaged_system} are formally established
in Proposition \ref{prop:Pass}. The proof is provided in Appendix
\ref{sec:Proof-of-Proposition}.
\begin{proposition}
\label{prop:Pass}Suppose Assumptions 1--3 hold with $\phi(t)=\xi(t)$,
and let 
$\mathcal{\varOmega}\left(\mathbf{0}_{d}\right)=\left\{ (\mathbf{x},\eta,\mathbf{z})|\left(\mathbf{1}_{N}\otimes I_{d}\right)^{\mathrm{T}}\mathbf{z}=\mathbf{0}_{d}\right\}.$
System \eqref{averaged_system} with initial values in $\mathcal{\varOmega}\left(\mathbf{0}_{d}\right)$
is globally uniformly asymptotically stable at the origin if \eqref{eq:vc_condition} holds and there exist
matrices $P_{2},P_{3}\in\mathbb{R}^{(N-1)\times(N-1)}$ with $P_{3}>0$
and positive scalars $\left\{ p_{11},p_{22},\delta\right\} $ such
that the following conditions hold:
\begin{equation}
\left[\begin{array}{cc}
p_{22}I_{N-1} & P_{2}\\
\star & P_{3}
\end{array}\right]>0,\quad\Phi_{11}<0,\quad\Phi_{2}<0,\label{eq:LMI_1}
\end{equation}
where $\Phi_{11} = 2\frac{\beta}{v}p_{11} - m p_{11} \alpha k + \delta M^{2}$,
$m = \min\{m_{i}\}$, $M = \max\{M_{i}\}$, and
\begin{equation}
\Phi_{2} =
\begin{bmatrix}
\Phi_{21} & \Phi_{22} & -\frac{1}{2}(p_{22}-p_{11})\alpha k I_{N-1} \\
\star & -\frac{1}{2}(P_{2} + P_{2}^{\mathrm{T}}) & -\frac{1}{2}\alpha k P_{2}^{\mathrm{T}} \\
\star & \star & -\delta I_{N-1}
\end{bmatrix}
\label{eq:phi_2}
\end{equation}
with
$
\Phi_{21} = -p_{22}\left(R^{\mathrm{T}} L R\right) - p_{22}\left(R^{\mathrm{T}} L^{\mathrm{T}} R\right)
 + \gamma P_{2}\left(R^{\mathrm{T}} L R\right) + \gamma\left(R^{\mathrm{T}} L^{\mathrm{T}} R\right) P_{2}^{\mathrm{T}},
\Phi_{22} = -p_{22} I_{N-1} + \gamma \left(R^{\mathrm{T}} L^{\mathrm{T}} R\right) P_{3}^{\mathrm{T}} - \left(R^{\mathrm{T}} L^{\mathrm{T}} R\right) P_{2},
$
in which $R\in\mathbb{R}^{N\times(N-1)}$ and $r=\frac{1}{\sqrt{N}}\mathbf{1}_{N}$ satisfy
\begin{equation}
R^{\mathrm{T}}r=\mathbf{0}_{N-1},\quad R^{\mathrm{T}}R=I_{N-1}.
\label{eq:rR_def}
\end{equation}
\end{proposition}

{\begin{remark}\label{rem:feasibility}
Consider the limiting case as $v\to\infty$, where $\xi(t)\equiv 1$ (see \eqref{eq:xi}) and 
the positive feedback term $\frac{\beta}{v}\xi^{-v}(t)\bar{\mathbf{x}}_{f}$ vanishes. 
Under condition \eqref{eq:vc_condition}, the time-varying derivatives in the averaged 
system \eqref{averaged_system} can be treated as vanishing perturbations that do not affect stability (see our proof in Appendix~\ref{sec:Proof-of-Proposition}), 
and the analysis reduces to the framework studied in \cite{kia2015distributed}. Notably, the Lyapunov function 
employed in Theorem 1 of \cite{kia2015distributed} corresponds to a special case of our candidate \eqref{Lyapunov} 
with the specific parameter selection: $p_{11}=\frac{1}{9}\alpha(\theta+1)$, $p_{22}=\alpha(\theta+1)$, $P_2=I_{N-1}$,
and $P_3=\frac{1}{\alpha}I_{N-1}$, where $\theta>0$ is a free tuning parameter. Although the conditions derived in Theorem 1 
of \cite{kia2015distributed} through completing the square are more concise, they are also {more conservative}.
 This observation suggests a practical parameter initialization and tuning procedure: (i) select initial parameters according 
 to Theorem 1 in \cite{kia2015distributed}, which guarantees LMI feasibility for sufficiently large $v$, and (ii) gradually 
 decrease $v$ from large values while fine-tuning other parameters to maintain the LMIs in Proposition \ref{prop:Pass} feasible.
\end{remark}}

Because the averaged system \eqref{averaged_system} is stable only
on an invariant set, the convergence of Algorithm \eqref{algo_extremum} is
established in the following senses.
\begin{definition}[SUAS and SPUAS w.r.t. $\mathcal{I}$]\cite{durr2013saddle} \label{def:spuas}
Consider a system $\dot{x}=F(x,t,\omega)$ parameterized by $\omega>0$, and let
$\mathcal{I}$ be a uniformly positively invariant set for the system. Its origin is
\emph{semi-globally uniformly asymptotically stable w.r.t. $\mathcal{I}$}
(SUAS w.r.t. $\mathcal{I}$) if for every compact
$\mathcal{K}\subseteq\mathcal{I}$ there is $\omega^{*}>0$ such that, for all
$\omega>\omega^{*}$ and $t_{0}\in\mathbb{R}$, every solution with
$x(t_{0})\in\mathcal{K}$ converges to the origin uniformly in $t_{0}$. It is
\emph{semi-globally practically uniformly asymptotically stable w.r.t.
$\mathcal{I}$} (SPUAS w.r.t. $\mathcal{I}$) when convergence holds only up to a
prescribed accuracy: for every $\varepsilon>0$ there is
$\omega^{*}(\varepsilon)>0$ such that the same solutions satisfy
$\limsup_{t\to\infty}\Vert x(t)\Vert\leq\varepsilon$.
\end{definition}

With this notion in place, we now establish the stability result
for Algorithm \eqref{algo_extremum}:
\begin{theorem}
\label{thm:1}Suppose Assumptions 1--3 hold with $\phi(t)=\xi(t)$,
and let $\mathcal{\varOmega}\left(\mathbf{0}_{d}\right)=\left\{ (\mathbf{x},\eta,\mathbf{z})|\left(\mathbf{1}_{N}\otimes I_{d}\right)^{\mathrm{T}}\mathbf{z}=\mathbf{0}_{d}\right\} $.
If conditions \eqref{eq:vc_condition} and \eqref{eq:LMI_1} are satisfied, then
algorithm \eqref{algo_extremum} renders the tracking error $\mathbf{x}-\boldsymbol{1}_{N}\otimes x^{*}$ SUAS w.r.t. $\mathcal{\varOmega}\left(\mathbf{0}_{d}\right)$, with $x_{i}\rightarrow x^{*}$ ($i=1,\dots,N$) at the rate of $\xi^{-1}(t)$.
\end{theorem}
\begin{proof}
By Proposition \ref{prop:Pass}, the origin of the
averaged system \eqref{averaged_system} is globally uniformly asymptotically
stable with respect to $\mathcal{\varOmega}\left(\mathbf{0}_{d}\right)$.
Furthermore, Lemma \ref{lem:bound} explicitly verifies Assumption A2 in \cite{durr2013lie},
while Assumptions A1, A3, and A4 hold straightforwardly for our system.
Therefore, Theorem 1 in \cite{durr2013lie} holds for
\eqref{eq:error_system} and its averaged system \eqref{averaged_system}:
for sufficiently large probing frequency $\omega$ the trajectories of the original
system \eqref{eq:error_system} remain arbitrarily close to those of the averaged system
\eqref{averaged_system} over any finite time interval.
Subsequently, following the same lines as the proof of Lemma 1 in
\cite{durr2013saddle}, we can prove that the origin of \eqref{eq:error_system}
is SPUAS w.r.t. $\mathcal{\varOmega}\left(\mathbf{0}_{d}\right)$.
Recalling the transformations
in \eqref{trans} we can provide bound on the convergence error as
\begin{equation*}
    \left\Vert \mathbf{x}-\boldsymbol{1}_{N}\otimes x^{*}\right\Vert \leq\xi^{-1}(t)\left\Vert \tilde{\mathbf{x}}_{f}\right\Vert .
\end{equation*}
Since $\xi^{-1}(t)\to0$, the tracking error is SUAS w.r.t. $\mathcal{\varOmega}\left(\mathbf{0}_{d}\right)$
and $x_{i}\rightarrow x^{*}$ ($i=1,\dots,N$) at the rate of $\xi^{-1}(t)$.
\end{proof}
\begin{remark}\label{rem:remark9}
It is important to note that system \eqref{averaged_system} is not \emph{globally}
asymptotically stable since its initial state must reside within the invariant
set $\mathcal{\varOmega}\left(\mathbf{0}_{d}\right)$. {To ensure the initial state lies in this invariant set, one can simply set $\mathbf{z}(0)=\mathbf{0}_{Nd}$. As $\mathbf{z}$ is internal to the algorithm, this restricts nothing in practice. Since the classical Lie bracket averaging theory in \cite{durr2013lie} requires the averaged system to be globally uniformly asymptotically stable, we invoke the results from \cite{durr2013saddle} to accommodate the analysis within an invariant set.}
\end{remark}
\begin{remark}\label{remark:c<0}
Although the condition $c<-3$ in Theorem \ref{thm:1} appears restrictive
at first glance, it constrains only Algorithm \eqref{algo_extremum}, the corresponding condition for Algorithm \eqref{algo_chirpy} is far weaker (see Remark \ref{remark:cp}). We emphasize two critical observations. First, for
time-invariant optimization problems (i.e., $\dot{\zeta}(t)=\ddot{\zeta}(t)=\dot{x}^{*}(t)=\ddot{x}^{*}(t)=0$), the algorithm
remains effective regardless of $c$, as the bounded dynamics assumption
reduces to triviality.  For instance, the time-invariant examples in Sections~\ref{subsec:Distributed-optimization-algorit} and~\ref{subsec:chirpy-sim}~1) trivially satisfy $c<-3$.
Second, and more fundamentally, Theorem \ref{thm:1}
establishes the analytical foundation for Section \ref{subsec:chirpy},
where we systematically extend this framework through time-varying
probing signals with parameterized growth rates.
\end{remark}

{
\begin{remark}\label{rem:single_vs_distributed}
In the single-agent setting, since $\nabla f(x^{*})=0$, the mean value theorem yields $\|\nabla f(x(t))\|=\|\nabla f(x(t))-\nabla f(x^{*})\|\leq\sigma_{2}\phi^{-1}(t)$ 
whenever $\|\tilde{x}(t)\|\leq\sigma_{1}\phi^{-1}(t)$. This inherent gradient decay property enables the uES framework to accommodate a broader class of convex functions than those considered in this paper, as well as exponential 
growth rates such as $\phi(t)=e^{\lambda t}$. In the distributed setting, however, $\nabla f_{i}(x^{*})\neq 0$ in general, 
and this decay property no longer holds. As a consequence, strong convexity of $\sum_{i}f_{i}$ becomes necessary, 
and a key requirement in our stability proof is that $\dot{\xi}(t)\to 0$ as $t\to\infty$, which is guaranteed by \eqref{eq:vc_condition} and serves to compensate for the absence of gradient decay. We note that this requirement precludes exponential growth functions such as $\phi(t)=e^{\lambda t}$.
\end{remark}
}

\section{Distributed Optimization Algorithms via uES with Chirpy Probing}
\label{subsec:chirpy}
Building on the asymptotic convergence guarantees of constant-frequency
probing, we now extend the framework to incorporate chirpy probing
signals with time-varying frequencies. This unified formulation enables
flexible convergence rate customization---achieving asymptotic, exponential,
or prescribed-time convergence through strategic selection of signal
growth parameters.  Intuitively, a faster-growing $\phi(t)$ acts as an accelerating clock, so the probe has to accelerate too: under the time-scale transformation the chirped dither returns to a constant frequency, and the constant-frequency Lie-bracket analysis still applies.
Inspired by \cite{yilmaz2024asymptotic}, we now
introduce three types of time-varying functions and parameters as
shown in Table \ref{table_chirpy}, in which the time-varying functions $\phi(t)$ and parameter $p$ are consistent with Table II in \cite{yilmaz2024asymptotic}, but the convergence conditions in the last column are different. With these, the algorithm \eqref{algo_extremum}
is modified as
\vspace{-0.5em}
\begin{equation}
\begin{aligned}
\dot{\mathbf{x}} &= \phi^{p}(t)\sum_{s=1}^{d}\sqrt{\alpha\omega_{s}}
    \cos\left( \omega_{s} \cdot \left( t_{0} + \rho \left( \phi^{q}(t) - 1 \right) \right) \right. \\
    & \qquad \left. + k\phi(t) \left( h(\mathbf{x},\zeta(t)) - \eta \right) \right) \otimes e_{s} \\
    &\quad + \left( -\phi^{p+1}(t)\mathbf{L}\mathbf{x} - \phi^{p}(t)\mathbf{z} \right), \\
\dot{\eta} &= \left( -\omega_{h}\eta + \omega_{h}h(\mathbf{x},\zeta(t)) \right) \phi^{p+1}(t), \\
\dot{\mathbf{z}} &= \gamma\phi^{p+2}(t)\mathbf{L}\mathbf{x}.
\end{aligned}
\label{algo_chirpy}
\end{equation}
\vspace{-0.3em}
The stability of Algorithm \eqref{algo_chirpy} is formalized in Theorem
\ref{thm:chirpy}.
\begin{table}
\centering
\setlength{\tabcolsep}{3pt} 
\begin{tabular}{|c|c|c|c|}
\hline
 & $\phi(t)$  & $p$  & Conditions \\
\hline
\makecell[c]{Asymptotic \\ ES}  & $\left(1+\beta\left(t-t_{0}\right)\right)^{\frac{1}{v}}$  & $q-v-1$  & $\begin{aligned} & \beta,v>0 \\ & q\geq2 \\ & \rho=v/(\beta q) \end{aligned}$ \\
\hline
\makecell[c]{Exponential \\ ES}  & $e^{\lambda\left(t-t_{0}\right)}$  & $q-1$  & $\begin{aligned} & \lambda>0 \\ & q\geq2 \\ & \rho=1/(\lambda q) \end{aligned}$ \\
\hline
\makecell[c]{Prescribed-time \\ ES}  & {$\left(\frac{T}{T+t_{0}-t}\right)^{\frac{1}{\varrho}}$}  & $q+\varrho-1$  & $\begin{aligned} & T,\varrho>0 \\ & q\geq2 \\ & \rho=\varrho T/q \end{aligned}$ \\
\hline
\end{tabular}
\caption{Time-varying functions used in \eqref{algo_chirpy}}
\label{table_chirpy}
\end{table}
\begin{theorem}
\label{thm:chirpy}Suppose Assumptions 1--3 hold, and let $\mathcal{\varOmega}\left(\mathbf{0}_{d}\right)=\left\{ (\mathbf{x},\eta,\mathbf{z})|\left(\mathbf{1}_{N}\otimes I_{d}\right)^{\mathrm{T}}\mathbf{z}=\mathbf{0}_{d}\right\} $.
Then algorithm \eqref{algo_chirpy} renders the tracking error $\mathbf{x}-\boldsymbol{1}_{N}\otimes x^{*}$ SUAS w.r.t. $\mathcal{\varOmega}\left(\mathbf{0}_{d}\right)$, with $x_{i}\rightarrow x^{*}$ ($i=1,\dots,N$) asymptotically, exponentially, or in prescribed time $T$ at the rate of $\phi^{-1}(t)$ given in Table \ref{table_chirpy},
if the LMIs in \eqref{eq:LMI_1} ($\beta/v$
is replaced by $1/(q\rho)$) hold together with

\begin{equation}\label{eq:qc_condition}
q\geq2,\qquad c-p<-2.
\end{equation}
If $\phi$ is saturated at $\bar{\phi}$ before $t_{0}+T$, the rate stops at $\bar{\phi}^{-1}$: an $O(\bar{\phi}^{-1})$ neighbourhood of $x^{*}$ is reached within the prescribed time and maintained.
\end{theorem}
\begin{proof}
We apply transformations analogous to \eqref{trans} by replacing the scaling function $\xi(t)$ with $\phi(t)$
 which convert \eqref{algo_chirpy} to
\allowdisplaybreaks[4]
\begin{equation}
\begin{aligned}
\dot{\tilde{\mathbf{x}}}_{f} &= \frac{\dot{\phi}(t)}{\phi(t)}\tilde{\mathbf{x}}_{f}
    + \phi^{p+1}(t)\sum_{s=1}^{d}\sqrt{\alpha\omega_{s}} \\
&\quad \times \cos\biggl( \omega_{s} \cdot \left( t_{0} + \rho \left( \phi^{q}(t) - 1 \right) \right) \\
& \qquad + k\phi(t)h_{f}(\tilde{\mathbf{x}}_{f},\zeta(t)) - k\tilde{\eta}_{f} \biggr) \otimes e_{s} \\
&\quad + \Bigl( -\phi^{p+1}(t)\mathbf{L}\tilde{\mathbf{x}}_{f} - \phi^{p+1}(t)\tilde{\mathbf{z}}_{f} - \phi^{p+1}(t)\mathbf{z}^{*}(t) \Bigr) \\
&\quad - \phi(t)\boldsymbol{1}_{N}\otimes\dot{x}^{*}(t) \\
\dot{\tilde{\eta}}_{f} &= \left(\frac{\dot{\phi}(t)}{\phi(t)} - \omega_{h}\phi^{p+1}(t)\right)\tilde{\eta}_{f}  + \omega_{h}h_{f}(\tilde{\mathbf{x}}_{f},\zeta(t))\phi^{p+2}(t) \\
&\quad - \phi(t)\dot{h}^{*}(t) \\
\dot{\tilde{\mathbf{z}}}_{f} &= \gamma\phi^{p+1}(t)\mathbf{L}\tilde{\mathbf{x}}_{f} - \dot{\mathbf{z}}^{*}(t)
\end{aligned}
\label{sys_trans_chirpy}
\end{equation}
where $h_{f}(\tilde{\mathbf{x}}_{f})$ is defined by \eqref{eq:h_f} and $h^{*}(t)=h(\boldsymbol{1}_{N}\otimes x^{*}(t),\zeta(t))$.
We now introduce the following time dilation and contraction transformations:
\begin{equation*}
\begin{aligned}
&\tau=t_{0}+\rho\left(\phi^{q}(t)-1\right),\\
&t=\phi^{\mathrm{inv}}\left(\left(1+\left(\tau-t_{0}\right)/\rho\right)^{\frac{1}{q}}\right),
\end{aligned}
\end{equation*}
where $\phi^{\mathrm{inv}}(\cdot)$ denotes the inverse of the function $\phi(\cdot)$.
Each $\phi(t)$ function in Table \ref{table_chirpy} with the corresponding
the parameters $p,q$ yields
\[
\begin{aligned}
\dot{\phi}(t) &= \phi^{p-q+2}(t)/(\rho q),\\
\frac{\mathrm{d}\tau}{\mathrm{d}t} &= \phi^{p+1}(t).
\end{aligned}
\]
For any $y(t)$ we write its $\tau$-domain counterpart as
$y_{\tau}(\tau)\triangleq y\left(\phi^{\mathrm{inv}}\left(\left(1+\left(\tau-t_{0}\right)/\rho\right)^{\frac{1}{q}}\right)\right)$;
in particular $\phi_{\tau}(\tau)=\left(1+\left(\tau-t_{0}\right)/\rho\right)^{\frac{1}{q}}$.
Through these transformations, we rewrite \eqref{sys_trans_chirpy}
in dilated $\tau$-domain as
\begin{align}
&\frac{\mathrm{d}\tilde{\mathbf{x}}_{f}}{\mathrm{d}\tau}=  \frac{1}{\rho q}\phi_{\tau}^{-q}(\tau)\tilde{\mathbf{x}}_{f}+\left(-\mathbf{L}\tilde{\mathbf{x}}_{f}-\tilde{\mathbf{z}}_{f}-\mathbf{z}_{\tau}^{*}(\tau)\right)\nonumber\\
&\;-\phi_{\tau}(\tau)\boldsymbol{1}_{N}\otimes\dot{x}_{\tau}^{*}(\tau)\nonumber\\
&\;+\sum_{s=1}^{d}\sqrt{\alpha\omega_{s}}\cos\left(\omega_{s}\tau+k\phi_{\tau}(\tau)h_{f}(\tilde{\mathbf{x}}_{f},\zeta_{\tau}(\tau))-k\tilde{\eta}_{f}\right)\otimes e_{s},\nonumber\\
&\frac{\mathrm{d}\tilde{\eta}_{f}}{\mathrm{d}\tau}=  \left(\frac{1}{\rho q}\phi_{\tau}^{-q}(\tau)-\omega_{h}\right)\tilde{\eta}_{f}+\omega_{h}h_{f}(\tilde{\mathbf{x}}_{f},\zeta_{\tau}(\tau))\phi_{\tau}(\tau)\nonumber\\
&\;-\phi_{\tau}(\tau)\dot{h}_{\tau}^{*}(\tau),\nonumber\\
&\frac{\mathrm{d}\tilde{\mathbf{z}}_{f}}{\mathrm{d}\tau}=  \gamma\mathbf{L}\tilde{\mathbf{x}}_{f}-\dot{\mathbf{z}}_{\tau}^{*}(\tau).
\label{sys_chripy_tau}
\end{align}
The growth bounds in Assumption 2 are rewritten as:
\[
\begin{aligned}
 & \left\|\frac{\mathrm{d}}{\mathrm{~d} \tau} \zeta_\tau(\tau)\right\|+\phi_\tau^{p+1}(\tau)\left\|\frac{\mathrm{d}^2}{\mathrm{~d} \tau^2} \zeta_\tau(\tau)\right\| \\
+{} & \left\|\frac{\mathrm{d}}{\mathrm{~d} \tau} x_\tau^*(\tau)\right\|+\phi_\tau^{p+1}(\tau)\left\|\frac{\mathrm{d}^2}{\mathrm{~d} \tau^2} x_\tau^*(\tau)\right\| \leqslant M_2 \phi_\tau^{c-p-1}(\tau).
\end{aligned}
\]
Note that system \eqref{sys_chripy_tau} has the similar form as system
\eqref{eq:error_system} with $q=v$, $\rho=1/\beta$, so following
the same arguments in the proof of Theorem \ref{thm:1} we can prove that the origin of system \eqref{sys_trans_chirpy} is SPUAS w.r.t. $\mathcal{\varOmega}\left(\mathbf{0}_{d}\right)$.
Finally, recalling the transformation $\tilde{\mathbf{x}}_{f} = \phi(t)(\mathbf{x} - \mathbf{1}_N \otimes x^*(t))$  we can
conclude that the tracking error is SUAS w.r.t. $\mathcal{\varOmega}\left(\mathbf{0}_{d}\right)$
and $x_{i}\rightarrow x^{*}$ ($i=1,\dots,N$) at the rate of $\phi^{-1}(t)$.
Let $\bar{t}$ be the instant at which $\phi$ reaches $\bar{\phi}$. Past $\bar{t}$ the dilated
time is continued as $\tau=\tau(\bar{t})+\bar{\phi}^{p+1}(t-\bar{t})$, which extends
$\mathrm{d}\tau/\mathrm{d}t=\phi^{p+1}(t)$ without a jump; freezing $\tau$ as well would arrest the
probing in \eqref{algo_chirpy}. Since $\dot{\phi}=0$, the
destabilising term $\frac{1}{\rho q}\phi_{\tau}^{-q}(\tau)\tilde{\mathbf{x}}_{f}$
of \eqref{sys_chripy_tau} disappears, and with it the $1/(q\rho)$ appearing in $\Phi_{11}$ of \eqref{eq:LMI_1}.
The LMIs are therefore only relaxed, $\tilde{\mathbf{x}}_{f}$ stays bounded, and
$\Vert\mathbf{x}-\boldsymbol{1}_{N}\otimes x^{*}\Vert=\Vert\tilde{\mathbf{x}}_{f}\Vert/\bar{\phi}$.
\end{proof}
\begin{remark}\label{remark:cp}
Compared to the condition $c<-3$ in Algorithm \eqref{algo_extremum},
in Algorithm \eqref{algo_chirpy} we only need \eqref{eq:qc_condition}, where $p$
is a tunable parameter as shown in Table \ref{table_chirpy}. By selecting sufficiently large $p$, this condition accommodates any bounded $c$. The theory thus imposes no fixed upper bound on the admissible rate of variation.
\end{remark}

\section{Numerical simulation}

\label{subsec:Numerical-Simulation}

This section presents comprehensive numerical validation of the proposed
algorithms in \eqref{algo_extremum} and \eqref{algo_chirpy}.\footnote{The code reproducing all numerical results in this section, the LMI verification programs, and further examples, including comparisons with the literature and a smart-grid application, are available at \url{https://github.com/Xuebin-Li-hit/Time-varyingES}.} We
implement a 5-agent network with directed communication topology shown
in Figure \ref{topology}, maintaining strong connectivity throughout
simulations.

\subsection{Distributed optimization algorithms via uES with constant-frequency
probing}

\label{subsec:Distributed-optimization-algorit}

\begin{figure*}[htbp]
\centering
\begin{subfigure}[c]{0.30\textwidth}
\centering \includegraphics[width=\linewidth]{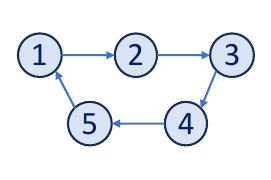}
\caption{}
\label{topology}
\end{subfigure}
\hfill
\begin{subfigure}[c]{0.33\textwidth}
\centering \includegraphics[width=\linewidth]{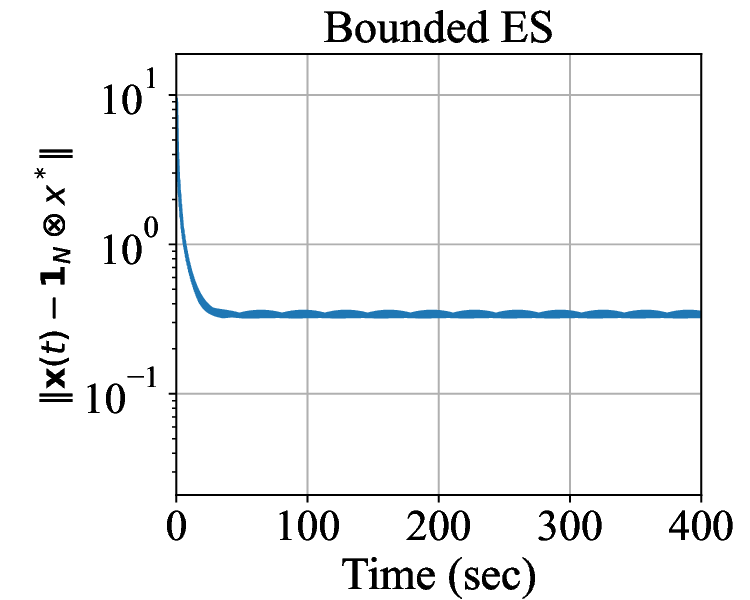}
\caption{}
\label{beta:subfig1}
\end{subfigure}
\hfill
\begin{subfigure}[c]{0.33\textwidth}
\centering \includegraphics[width=\linewidth]{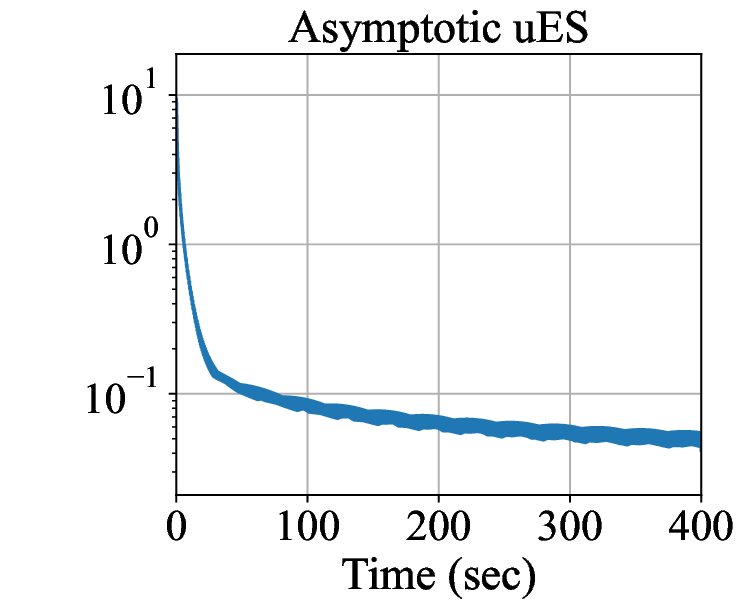}
\caption{}
\label{beta:subfig2}
\end{subfigure}
\caption(a) Directed communication graph for the 5-agent network. Tracking
error $\|\mathbf{x}(t)-\boldsymbol{1}_{N}\otimes x^{*}\|$ on a semilogarithmic scale
($d=3$) for distributed optimization via ES with constant-frequency probing~(\ref{algo_extremum}):
(b) $\beta=0$ configuration, (c) $\beta=0.2$ configuration.
\label{beta}
\end{figure*}

We first evaluate Algorithm \eqref{algo_extremum} in the
multi-dimensional setting $d=3$ using the following cost functions:
\begin{equation*}
f_{i}(x)=\|x-c_{i}\|^{2}+\ln\!\big(1+\|x-c_{i}\|^{2}\big),\quad i=1,2,\dots,5,
\end{equation*}
where the minimizers $c_{1},\dots,c_{5}$ and the initial states
$x_{1}(0),\dots,x_{5}(0)$ are, respectively, the columns of
\[
\left[\begin{array}{rrrrr}
1 & 2 & 3 & 4 & 5\\
2 & 1 & 3 & 2 & 4\\
1 & 3 & 2 & 4 & 3
\end{array}\right],\qquad
\left[\begin{array}{rrrrr}
-1 & 0 & 1 & 4 & 5\\
0 & 3 & -1 & 5 & 1\\
2 & -1 & 4 & 0 & 5
\end{array}\right],
\]
so that $f(x)=\sum_{i=1}^{5}f_{i}(x)$ is minimized at
$x^{*}\approx(3,2.4124,2.5876)^{\top}$.
The parameters are $v=2$, $k=1$, $\gamma=0.1$, $\alpha=0.4$, $\omega_{h}=10$,
$\omega=10$ and $(\hat{\omega}_{1},\hat{\omega}_{2},\hat{\omega}_{3})=(3,5,7)$. For these values the LMIs in \eqref{eq:LMI_1} are feasible, a certificate being $p_{11}=2.33$, $p_{22}=1$, $\delta=0.072$ together with matrices $P_{2},P_{3}$ returned by a standard semidefinite program.
To establish a comparative baseline, we implement the bounded ES algorithm by configuring $\beta=0$, and implement the asymptotic
uES algorithm with enhanced configuration $\beta=0.2$. Note that $\beta=0$ gives $\xi(t)\equiv1$, whereby \eqref{algo_extremum} reduces to a bounded ES scheme \cite{scheinker2014extremum} with PI consensus, which attains only practical (biased) convergence. As Figure \ref{beta} shows, with $\beta=0$ the error stalls at
$3.3\times10^{-1}$, the residual bias inherent to a constant probing
amplitude, whereas with $\beta=0.2$ it decreases monotonically to
$4.2\times10^{-2}$ at $t=400$ without settling on a floor, corroborating the
unbiased convergence of Theorem \ref{thm:1}.

\subsection{Distributed optimization algorithms via uES with chirpy probing}

\label{subsec:chirpy-sim}

\begin{figure*}[htbp]
\centering
\begin{subfigure}[t]{0.32\textwidth}
\centering \includegraphics[width=\linewidth]{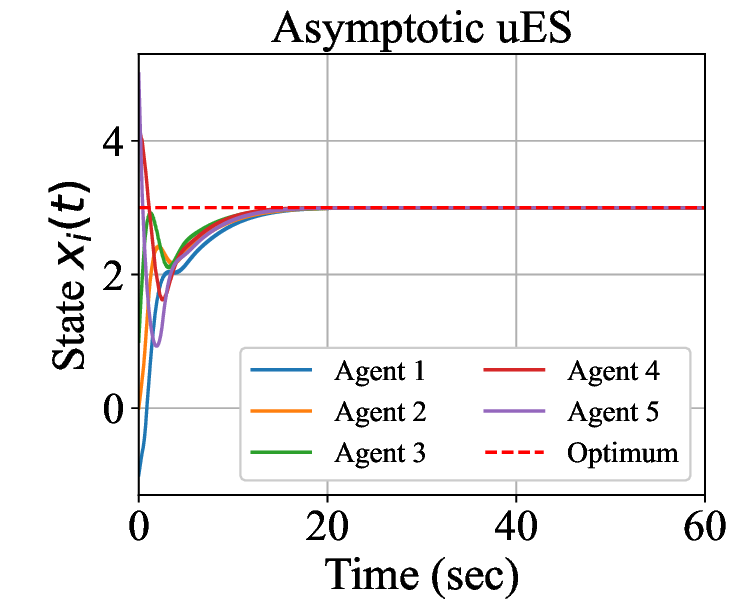}
\caption{}
\label{chirpy_invariant:subfig1}
\end{subfigure}
\hfill
\begin{subfigure}[t]{0.32\textwidth}
\centering \includegraphics[width=\linewidth]{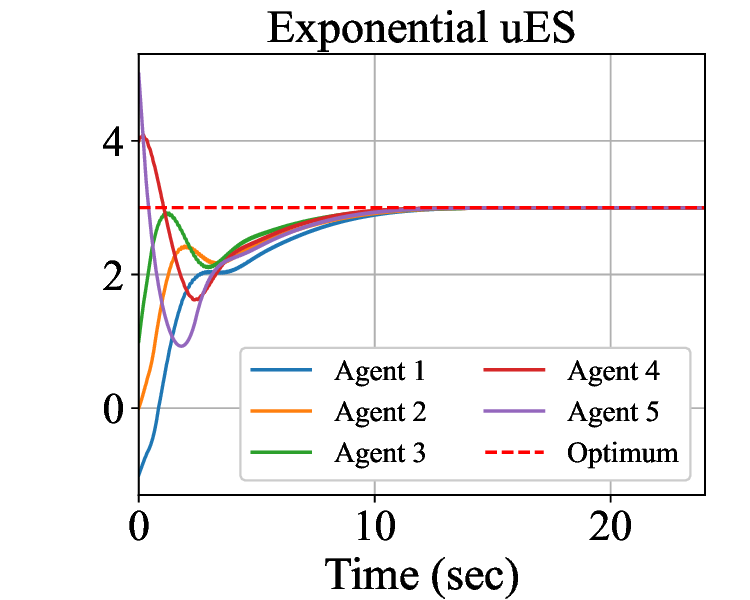}
\caption{}
\label{chirpy_invariant:subfig2}
\end{subfigure}
\hfill
\begin{subfigure}[t]{0.32\textwidth}
\centering \includegraphics[width=\linewidth]{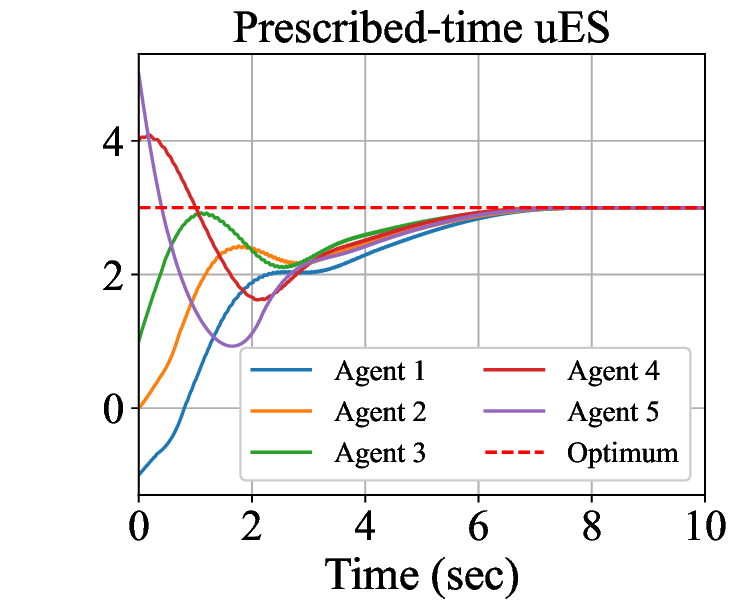}
\caption{}
\label{chirpy_invariant:subfig3}
\end{subfigure}

\caption{uES for time-invariant distributed optimization with chirpy probing.
(a)$\phi(t)=\left(1+\beta\left(t-t_{0}\right)\right)^{\frac{1}{v}}$
with $\beta=0.05,v=0.5$, (b) $\phi(t)=e^{\lambda\left(t-t_{0}\right)}$
with $\lambda=0.1$, (c) $\phi(t)=\left(\frac{T}{T+t_{0}-t}\right)^{\frac{1}{\varrho}}$
with $T=10,\varrho=1$.}
\label{chirpy_invariant}
\end{figure*}

\textbf{1) Time-invariant extrema:} We first consider the distributed
algorithm presented by \eqref{algo_chirpy} for time-invariant extrema.
We take $d=1$ with $f_{i}(x)=(x-i)^{2}+\ln(1+(x-i)^{2})$,
$i=1,\dots,5$, minimized at $x^{*}=3$, and $\mathbf{x}(0)=[-1,0,1,4,5]^{\top}$.
In the following simulations, parameters $q=2,k=32,\gamma=0.05,\alpha=0.0125,\omega=40$,
and $\omega_{h}=8$ remain constant. The three configurations share the
value $1/(q\rho)=0.1$, hence the same LMIs \eqref{eq:LMI_1}, which are feasible
for these gains, a certificate being $p_{11}=2.38$, $p_{22}=1$, $\delta=0.073$.
Figure \ref{chirpy_invariant} demonstrates the performance of three
chirpy probing configurations for time-invariant distributed optimization.
The left, middle, and right plots illustrate the convergence behavior
corresponding to the asymptotic, exponential, and prescribed-time
configurations, respectively.

\begin{figure}[htbp]
\centering \begin{subfigure}[t]{0.48\textwidth} \centering \includegraphics[width=1\columnwidth]{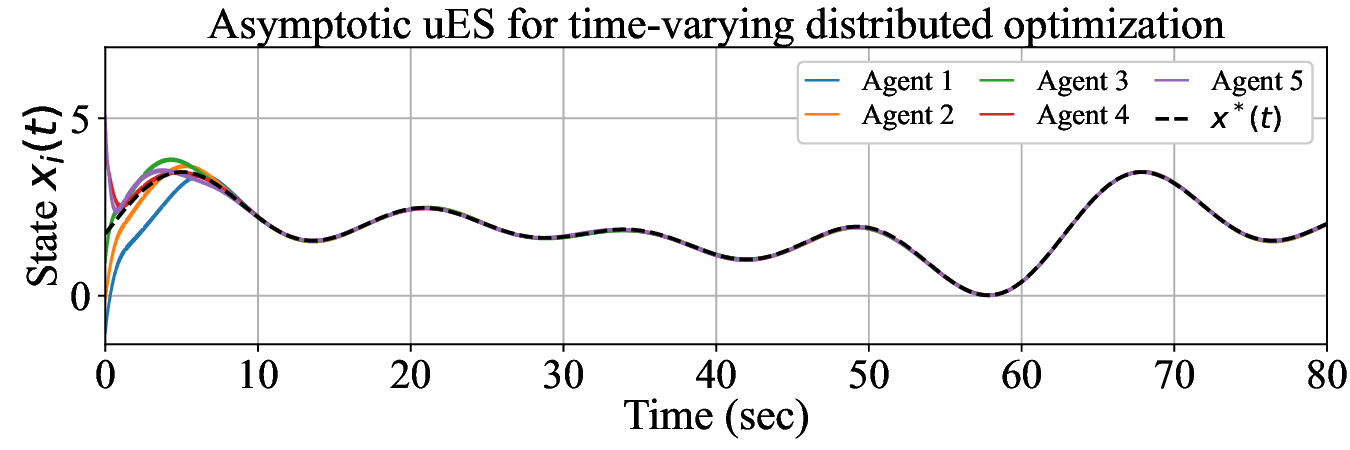}
\caption{}
\label{chirpy_varying:subfig1} \end{subfigure}

\begin{subfigure}[t]{0.48\textwidth} \centering \includegraphics[width=1\columnwidth]{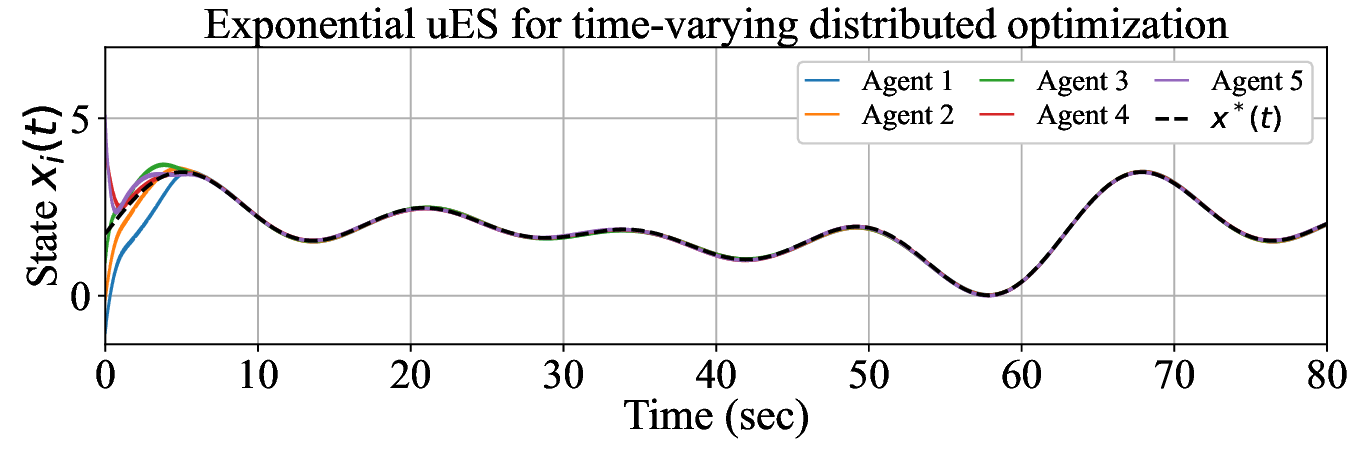}
\caption{}
\label{chirpy_varying:subfig2} \end{subfigure}

\begin{subfigure}[t]{0.48\textwidth} \centering \includegraphics[width=1\columnwidth]{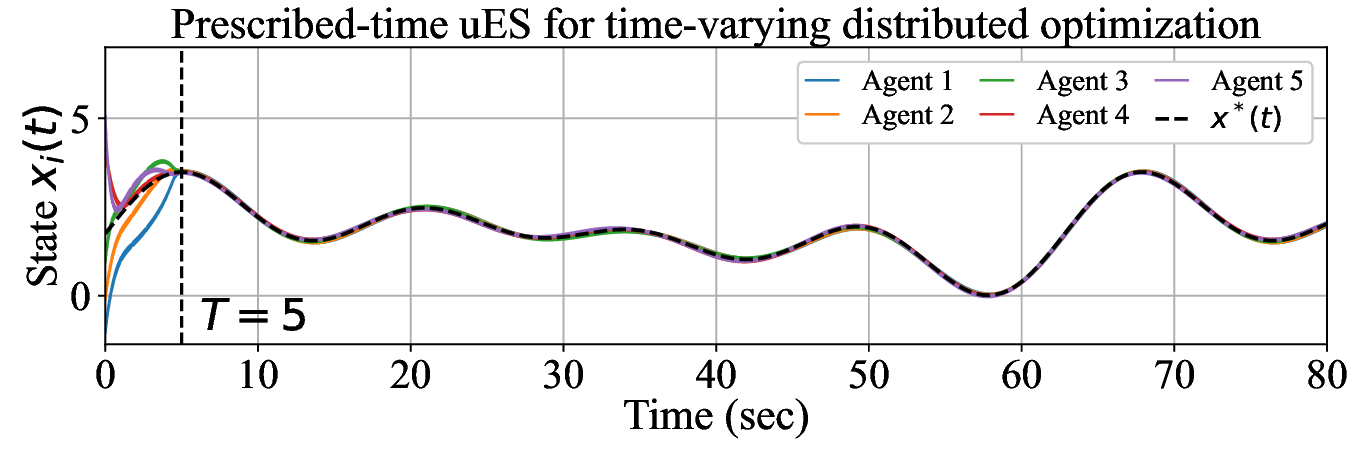}
\caption{}
\label{chirpy_varying:subfig3} \end{subfigure}

\caption{uES for time-varying distributed optimization with chirpy probing}
\label{chirpy_varying}
\end{figure}

\textbf{2) Time-varying extrema:} We next validate the proposed algorithm
\eqref{algo_chirpy} for time-varying optimization scenarios. The
time-varying local cost functions are defined as:
$f_{1} = \left(x - 0.75 - 0.6\sin(0.1t)\right)^{2}$, $f_{2} = \left(x - 1.25 - 1.8\sin(0.2t)\right)^{2}$, $f_{3} = \left(x - 1.75 - 3\sin(0.3t)\right)^{2}$, $f_{4} = \left(x - 2.25 - 2.4\sin(0.1t)\right)^{2}$, $f_{5} = \left(x - 2.75 - 3\sin(0.4t)\right)^{2}$.
The algorithm parameters are configured as: $k=1$, $\gamma=0.05$,
$\omega=400$, $\omega_{h}=8$, $\alpha=1.6$. Initial values
are set to $\mathbf{x}(0)=[-1,0,1,4,5]^{\top}$, $\eta(0)=\mathbf{z}(0)=\boldsymbol{0}$. The growth functions are those of Table \ref{table_chirpy}, with
$\beta=0.1$, $v=0.5$, $q=4$ for the asymptotic configuration, $\lambda=0.2$,
$q=4$ for the exponential one and $T=5$, $\varrho=2$, $q=2$ for the
prescribed-time one, so that $p>2$ in all three, as required by $c-p<-2$ with
$c=0$. The LMIs \eqref{eq:LMI_1} are again feasible, a certificate being
$p_{11}=2.51$, $p_{22}=1$, $\delta=1.70$. The tracking performances are illustrated
in Figure \ref{chirpy_varying}.
Note that for the prescribed-time algorithm, the growth function $\phi(t)$
is held constant beyond $t=4.5$. Specifically: $\phi(t)=\phi(4.5)$
for $t>4.5$ (see Theorem \ref{thm:chirpy}), as is likewise required
in the single-agent designs of \cite{yilmaz2024asymptotic}.
\section{Conclusion}
\label{sec:Conclusion}
This work established a distributed optimization framework that enables
gradient-free tracking of time-varying extrema through uES principles,
eliminating the need for explicit mathematical formulations of objective
functions and relying solely on real-time measurements of function outputs.
Theoretical guarantees have been derived through Lie bracket averaging and
Lyapunov analysis, with stability conditions formulated as verifiable LMIs.
Numerical validations have confirmed the approach's effectiveness in both
static and dynamic scenarios, where the unbiased scheme removes
the residual error of its bounded-ES counterpart.  Sampled-data, delayed and quantized
implementation of the suggested algorithms, as studied e.g.~in~\cite{bo2024quantization,liu2018event}, may be a topic of future research.

\appendices

\section{Proof of Lemma \ref{lem:bound}}
\label{sec:Proof-of-Lemmabound}

To avoid complicated vector operator notation, we assume $d = l = 1$ when
proving the boundedness of all expressions except for the Lie bracket
$\left[b_{1s}(\varsigma,t),b_{2s}(\varsigma,t)\right]$. This simplifies the
derivation considerably; the boundedness of vector quantities is established by
verifying the boundedness of each of their components.

Defining
\[
h_{fi}(\tilde{x}_{fi},\zeta(t)) = f_{i}(\xi^{-1}(t)\tilde{x}_{fi}+x^{*}(t),\zeta(t))-f_{i}(x^{*}(t),\zeta(t)),
\]
we have $h_{f}(\tilde{\mathbf{x}}_{f},\zeta(t))=\operatorname{col}\left\{ h_{fi}(\tilde{x}_{fi},\zeta(t))\right\} _{i=1}^{N}$.
Since $x_{i}=\xi^{-1}(t)\tilde{x}_{fi}+x^{*}(t)$, item~4 of
Lemma \ref{lem:Lipschitz_grad} and item~4 of Lemma \ref{lem:strong_convex},
both applied at $x=x^{*}(t)$ and $y=x_{i}$, bracket $h_{fi}$ from both sides:
\[
a_{i}+\tfrac{m_{i}}{2}\Vert x_{i}-x^{*}\Vert^{2}\;\leq\;h_{fi}\;\leq\;a_{i}+\tfrac{M_{i}}{2}\Vert x_{i}-x^{*}\Vert^{2},
\]
where $a_{i}=\frac{\partial f_{i}}{\partial x_{i}}(x^{*}(t),\zeta(t))^{\mathrm{T}}(x_{i}-x^{*}(t))$.
The upper bound alone controls the signed scalar $h_{fi}$ but not its modulus;
the matching lower bound, supplied by the $m_{i}$-strong convexity of
Assumption 1, is what makes $\Vert h_{fi}\Vert$ accessible. As $0<m_{i}\leq M_{i}$,
both endpoints lie in $\left[a_{i},\,a_{i}+\tfrac{M_{i}}{2}\Vert x_{i}-x^{*}\Vert^{2}\right]$,
so $\Vert h_{fi}\Vert\leq\vert a_{i}\vert+\tfrac{M_{i}}{2}\Vert x_{i}-x^{*}\Vert^{2}$.
Bounding $\vert a_{i}\vert$ by Cauchy--Schwarz and substituting
$x_{i}-x^{*}(t)=\xi^{-1}(t)\tilde{x}_{fi}$ then yields
\begin{equation}
\begin{array}{@{}r@{}l@{}}
\xi(t)\Vert h_{fi}\Vert & {}\leq \left\Vert \frac{\partial f_{i}}{\partial x_{i}}(x^{*}(t),\zeta(t)) \right\Vert \|\tilde{x}_{fi}\| \\
& \quad + \frac{M_{i}}{2}\xi^{-1}(t)\|\tilde{x}_{fi}\|^{2}, \quad i=1,\dots,N.
\end{array}
\label{eq:lip_grad}
\end{equation}
Here $\tilde{x}_{fi}$ ranges over the compact set $\mathcal{X}$ and
$\xi^{-1}(t)\leq1$, so the right-hand side of \eqref{eq:lip_grad} is bounded as
soon as $\|\frac{\partial f_{i}}{\partial x_{i}}(x^{*}(t),\zeta(t))\|$ is;
the latter holds because $\|x^{*}(t)\|$ and $\|\zeta(t)\|$ are bounded by
Assumption 2 and $f_{i}$ is $C^{3}$ by Assumption 1. Hence
$\xi(t)\Vert h_{fi}\Vert$ is bounded, which is what the $\tilde{\eta}_{f}$
block of $b_{0}$ requires. Moreover
\[
\|x_{i}\|\leq\|x_{i}-x^{*}(t)\|+\|x^{*}(t)\|\leq\xi^{-1}(t)\|\tilde{x}_{fi}\|+\|x^{*}(t)\|,
\]
so $\|x_{i}\|$ is bounded as well, and with the smoothness in Assumption 1 the
partial derivatives of $f_{i}(x_{i},\zeta(t))$ with respect to $x_{i}$ and
$\zeta$ up to the third order are all bounded. In particular
$\mathbf{z}^{*}(t)$ is bounded. Since
\[
\frac{\mathrm{d}}{\mathrm{d}t}\left(\frac{\partial f_{i}(x^{*}(t),\zeta(t))}{\partial x_{i}}\right)
= \frac{\partial^{2}f_{i}}{\partial x_{i}^{2}} \dot{x}^{*}(t) + \frac{\partial^{2}f_{i}}{\partial x_{i}\partial\zeta} \dot{\zeta}(t),
\]
Assumption 2 yields a $k_{1}>0$ such that
\begin{equation}\label{bound_z}
    \left\Vert \frac{\mathrm{d}}{\mathrm{d}t}\left(\frac{\partial G}{\partial\mathbf{x}}(\boldsymbol{1}_{N}\otimes x^{*}(t),\zeta(t))\right)\right\Vert \leq k_{1}\xi^{c}(t).
\end{equation}
Likewise, from
\[
\frac{\mathrm{d}}{\mathrm{d}t}\left(f_{i}(x^{*}(t),\zeta(t))\right)
= \frac{\partial f_{i}}{\partial x_{i}} \dot{x}^{*}(t) + \frac{\partial f_{i}}{\partial\zeta} \dot{\zeta}(t)
\]
and Assumption 2 with $c\leq-1$, the term
$\xi(t)\dot{h}(\boldsymbol{1}_{N}\otimes x^{*}(t),\zeta(t))$ appearing in
$b_{0}$ is bounded. Collecting these facts,
$\|b_{0}(\varsigma,t)\|$, $\|b_{1s}(\varsigma,t)\|$ and
$\|b_{2s}(\varsigma,t)\|$ are uniformly bounded for $\varsigma\in\mathcal{X}$.

\emph{Spatial derivatives.} Noting $\frac{\partial x_{i}}{\partial\tilde{x}_{fi}}=\xi^{-1}(t)I_{d}$,
\begin{equation}
\begin{array}{@{}r@{}l@{}}
\left\Vert \frac{\partial\left(\xi(t)h_{fi}\right)}{\partial\tilde{x}_{fi}}\right\Vert & {}= \left\Vert \frac{\partial f_{i}}{\partial x_{i}}(x_{i},\zeta(t))\right\Vert \\
& {}= \left\Vert \frac{\partial f_{i}}{\partial x_{i}}(x_{i},\zeta(t)) - \frac{\partial f_{i}}{\partial x_{i}}(x_i^{*}(t),\zeta(t))\right\Vert \\
& {}\leq M_{i}\|x_{i}-x_i^{*}(t)\| \\
& {}\leq M_{i}\xi^{-1}(t)\|\tilde{x}_{fi}\|+M_{i}\|x^{*}(t)-x_i^{*}(t)\|,
\end{array}
\label{eq:partial_x}
\end{equation}
where the second equality uses $\frac{\partial f_{i}}{\partial x_{i}}(x_i^{*}(t),\zeta(t))=0$.
Hence $\|\frac{\partial b_{0}}{\partial\varsigma}\|$,
$\|\frac{\partial b_{1s}}{\partial\varsigma}\|$ and
$\|\frac{\partial b_{2s}}{\partial\varsigma}\|$ are bounded for
$\varsigma\in\mathcal{X}$.

\emph{Time derivatives.} By direct computation,
\[
\begin{array}{@{}r@{}l@{}}
 & \frac{\partial\left(\xi(t)h_{fi}\right)}{\partial t} = \frac{\beta}{v}\xi^{1-v}(t)h_{fi} + \xi(t)\Big(-\frac{\mathrm{d}}{\mathrm{d}t}f_{i}(x^{*}(t),\zeta(t)) \\
 & \quad + \frac{\partial f_{i}}{\partial x_{i}}(x_{i},\zeta(t))\Big(-\frac{\beta}{v}\xi^{-1-v}(t)\tilde{x}_{fi}+\dot{x}^{*}(t)\Big) \\
 & \quad + \frac{\partial f_{i}}{\partial\zeta}(x_{i},\zeta(t))\dot{\zeta}(t)\Big) \\
= {} & \frac{\beta}{v}\xi^{1-v}(t)h_{fi} \\
 & + \xi(t)\left(\frac{\partial f_{i}}{\partial x_{i}}(x_{i},\zeta(t)) - \frac{\partial f_{i}}{\partial x_{i}}(x^{*}(t),\zeta(t))\right)\dot{x}^{*}(t) \\
 & - \frac{\beta}{v}\xi^{-v}(t)\frac{\partial f_{i}}{\partial x_{i}}(x_{i},\zeta(t))\tilde{x}_{fi} \\
 & + \xi(t)\left(\frac{\partial f_{i}}{\partial\zeta}(x_{i},\zeta(t)) - \frac{\partial f_{i}}{\partial\zeta}(x^{*}(t),\zeta(t))\right)\dot{\zeta}(t).
\end{array}
\]
The second form makes the boundedness transparent: every term
carrying a positive power of $\xi(t)$ multiplies a difference evaluated at
$x_{i}$ and $x^{*}(t)$, and since $x_{i}-x^{*}(t)=\xi^{-1}(t)\tilde{x}_{fi}$
that difference is $O(\xi^{-1}(t)\|\tilde{x}_{fi}\|)$ by the $M_{i}$-Lipschitz
continuity of the derivatives of $f_{i}$. The growth of $\xi(t)$ therefore
cancels and no unbounded term survives. Term by term, with
$\tilde{x}_{fi}\in\mathcal{X}$ and
$\|\dot{x}^{*}\|,\|\dot{\zeta}\|\leq M_{2}\xi^{c}(t)$. The first term is
$O(\xi^{-v}(t))$ by \eqref{eq:lip_grad} and vanishes since $v\geq2$. In the
second, the $M_{i}$-Lipschitz continuity of $\frac{\partial f_{i}}{\partial x_{i}}$
bounds the bracket by $M_{i}\|x_{i}-x^{*}(t)\|=M_{i}\xi^{-1}(t)\|\tilde{x}_{fi}\|$,
so the term is at most $M_{i}M_{2}\xi^{c}(t)\|\tilde{x}_{fi}\|$. The third is
bounded by $\xi^{-v}(t)\leq1$ and the boundedness of the remaining factors. In the
fourth, $\|\frac{\partial^{2}f_{i}}{\partial x_{i}\partial\zeta}\|\leq M_{3}$,
which exists because $f_{i}$ is $C^{3}$ and $x_{i},\zeta$ stay bounded, bounds the
bracket by $M_{3}\xi^{-1}(t)\|\tilde{x}_{fi}\|$, leaving at most
$M_{3}M_{2}\xi^{c}(t)\|\tilde{x}_{fi}\|$. Since $\xi(t)\geq1$ and $c<0$, all
four terms are bounded.

The time derivatives of the other terms entering $b_{0}$ are
\[
\begin{array}{@{}r@{}l@{}}
\frac{\partial\left(\xi(t)\dot{f_{i}}(x^{*}(t),\zeta(t))\right)}{\partial t} = {} & \frac{\beta}{v}\xi^{1-v}(t)\dot{f_{i}}(x^{*}(t),\zeta(t)) \\
& + \xi(t)\Big(\dot{x}^{*\mathrm{T}}(t)\frac{\partial^{2}f_{i}}{\partial x_{i}^{2}}\dot{x}^{*}(t) \\
& \quad + \frac{\partial^{2}f_{i}}{\partial\zeta^{2}}\dot{\zeta}^{2}(t) + 2\frac{\partial^{2}f_{i}}{\partial x_{i}\partial\zeta}\dot{x}^{*}(t)\dot{\zeta}(t) \\
& \quad + \frac{\partial f_{i}}{\partial x_{i}}\ddot{x}^{*}(t) + \frac{\partial f_{i}}{\partial\zeta}\ddot{\zeta}(t)\Big), \\
\frac{\partial\left(\xi(t)\dot{x}^{*}(t)\right)}{\partial t} = {} & \frac{\beta}{v}\xi^{1-v}(t)\dot{x}^{*}(t)+\xi(t)\ddot{x}^{*}(t), \\
\frac{\partial\left(\frac{\beta}{v}\xi^{-v}(t)\tilde{\eta}_{f}\right)}{\partial t} = {} & -\frac{\beta^{2}}{v}\xi^{-2v}(t)\tilde{\eta}_{f}.
\end{array}
\]
Assumption 2 bounds $\|\dot{x}^{*}\|,\|\ddot{x}^{*}\|,\|\dot{\zeta}\|,\|\ddot{\zeta}\|$
by $M_{2}\xi^{c}(t)$ with $c\leq-1$, which \eqref{eq:vc_condition} supplies.
Among the terms carrying an explicit factor $\xi(t)$, those linear in the
derivatives of $x^{*}$ and $\zeta$ are $O(\xi^{1+c}(t))=O(1)$, while those
quadratic in them are $O(\xi^{1+2c}(t))$ and hence smaller; the remaining terms
decay as $\xi^{1-v}(t)$ or faster. All of them are therefore bounded. Finally,
\[
\begin{array}{@{}r@{}l@{}}
\frac{\mathrm{d}^{2}}{\mathrm{d}t^{2}}\left(\frac{\partial f_{i}(x^{*}(t),\zeta(t))}{\partial x_{i}}\right) = {} & \frac{\partial^{3}f_{i}}{\partial x_{i}^{3}}\left(\dot{x}^{*}(t)\right)^{2} \\
& + 2\frac{\partial^{3}f_{i}}{\partial x_{i}^{2}\partial\zeta}\dot{x}^{*}(t)\dot{\zeta}(t) \\
& + \frac{\partial^{3}f_{i}}{\partial x_{i}\partial\zeta^{2}}\left(\dot{\zeta}(t)\right)^{2} \\
& + \frac{\partial^{2}f_{i}}{\partial x_{i}^{2}}\ddot{x}^{*}(t) + \frac{\partial^{2}f_{i}}{\partial x_{i}\partial\zeta}\ddot{\zeta}(t),
\end{array}
\]
which is bounded because the derivatives of $f_{i}$ up to third order are
bounded; hence, by \eqref{eq:z_star}, $\|\ddot{\mathbf{z}}^{*}(t)\|$ is bounded,
which is what differentiating the $\tilde{\mathbf{z}}_{f}$ block of $b_{0}$
requires. Consequently
$\|\frac{\partial b_{0}}{\partial t}\|$, $\|\frac{\partial b_{1s}}{\partial t}\|$
and $\|\frac{\partial b_{2s}}{\partial t}\|$ are bounded for
$\varsigma\in\mathcal{X}$.

\emph{Mixed partial derivatives.} From \eqref{eq:partial_x},
\begin{equation}
\begin{array}{@{}r@{}l@{}}
\frac{\partial^{2}\left(\xi(t)h_{fi}\right)}{\partial t\partial\tilde{x}_{fi}} & {}= \frac{\partial^{2}f_{i}(x_{i},\zeta(t))}{\partial x_{i}\partial t} \\
& {}= \frac{\partial^{2}f_{i}}{\partial x_{i}^{2}}\frac{\partial x_{i}}{\partial t} + \frac{\partial^{2}f_{i}}{\partial x_{i}\partial\zeta}\dot{\zeta}(t) \\
& {}= \frac{\partial^{2}f_{i}}{\partial x_{i}^{2}}\left(-\frac{\beta}{v}\xi^{-1-v}(t)\tilde{x}_{fi}+\dot{x}^{*}(t)\right)\\
& \quad + \frac{\partial^{2}f_{i}}{\partial x_{i}\partial\zeta}\dot{\zeta}(t).
\end{array}
\label{eq:mix_partial}
\end{equation}
Bounding each of the three resulting terms separately, and using
$\|\frac{\partial^{2}f_{i}}{\partial x_{i}^{2}}\|\leq M_{i}$ together with
$\|\dot{x}^{*}\|,\|\dot{\zeta}\|\leq M_{2}\xi^{c}(t)$, gives
\[
\begin{array}{@{}r@{}l@{}}
\left\Vert \frac{\partial^{2}\left(\xi(t)h_{fi}\right)}{\partial t\partial\tilde{x}_{fi}}\right\Vert & {}\leq \frac{\beta}{v}\xi^{-1-v}(t)M_{i}\|\tilde{x}_{fi}\| \\
& \quad + \left(M_{i}+M_{3}\right)M_{2}\xi^{c}(t),
\end{array}
\]
with $M_{3}$ as above. All
three terms are bounded, so $\|\frac{\partial^{2}b_{0}}{\partial\varsigma\partial t}\|$,
$\|\frac{\partial^{2}b_{1s}}{\partial\varsigma\partial t}\|$ and
$\|\frac{\partial^{2}b_{2s}}{\partial\varsigma\partial t}\|$ are bounded for
$\varsigma\in\mathcal{X}$.

\emph{The Lie bracket.} Write $\vartheta_{i}(\varsigma)=k\xi(t)h_{fi}(\tilde{x}_{fi},\zeta(t))-k\tilde{\eta}_{f,i}$
for the probing phase of agent $i$, so that
$b_{1s}=\left[\sqrt{\alpha}\operatorname{col}\{\cos\vartheta_{i}\}_{i=1}^{N}\otimes e_{s};\ \mathbf{0}_{N+Nd}\right]$
and $b_{2s}=\left[-\sqrt{\alpha}\operatorname{col}\{\sin\vartheta_{i}\}_{i=1}^{N}\otimes e_{s};\ \mathbf{0}_{N+Nd}\right]$.
Since $\frac{\partial\vartheta_{i}}{\partial\tilde{x}_{fi}}=k\xi(t)\frac{\partial h_{fi}}{\partial\tilde{x}_{fi}}=k\frac{\partial f_{i}}{\partial x_{i}}$
by \eqref{eq:partial_x}, and $\frac{\partial\vartheta_{i}}{\partial\tilde{\eta}_{f,i}}=-k$,
\[
\begin{array}{@{}r@{}l@{}}
\frac{\partial b_{1s}}{\partial\varsigma} = {} & \left[\begin{array}{ccc}\Theta_{1} & \Theta_{2} & 0\\ 0 & 0 & 0\\ 0 & 0 & 0\end{array}\right],\quad
\frac{\partial b_{2s}}{\partial\varsigma} = \left[\begin{array}{ccc}\Theta_{3} & \Theta_{4} & 0\\ 0 & 0 & 0\\ 0 & 0 & 0\end{array}\right],
\end{array}
\]
where
\[
\begin{array}{@{}r@{}l@{}}
\Theta_{1} & {}= -\sqrt{\alpha}k\operatorname{diag}\left\{ \sin\vartheta_{i}\left(\frac{\partial f_{i}}{\partial x_{i}}\right)^{\mathrm{T}}\right\} _{i=1}^{N}\otimes e_{s}, \\
\Theta_{2} & {}= \sqrt{\alpha}k\operatorname{diag}\left\{ \sin\vartheta_{i}\right\} _{i=1}^{N}\otimes e_{s}, \\
\Theta_{3} & {}= -\sqrt{\alpha}k\operatorname{diag}\left\{ \cos\vartheta_{i}\left(\frac{\partial f_{i}}{\partial x_{i}}\right)^{\mathrm{T}}\right\} _{i=1}^{N}\otimes e_{s}, \\
\Theta_{4} & {}= \sqrt{\alpha}k\operatorname{diag}\left\{ \cos\vartheta_{i}\right\} _{i=1}^{N}\otimes e_{s}.
\end{array}
\]
Because the $\tilde{\eta}_{f}$ and $\tilde{\mathbf{z}}_{f}$ blocks of $b_{1s}$
and $b_{2s}$ vanish, only $\Theta_{1}$ and $\Theta_{3}$ contribute, and
\[
\begin{array}{@{}r@{}l@{}}
\left[b_{1s},b_{2s}\right] & {}= \frac{\partial b_{2s}}{\partial\varsigma}b_{1s}-\frac{\partial b_{1s}}{\partial\varsigma}b_{2s} \\
& {}= \left[\begin{array}{c} -\alpha k\operatorname{col}\left\{ \left(\cos^{2}\vartheta_{i}+\sin^{2}\vartheta_{i}\right)\frac{\partial f_{i}}{\partial x_{is}}\right\} _{i=1}^{N}\otimes e_{s}\\ \mathbf{0}_{N+Nd}\end{array}\right] \\
& {}= \left[\begin{array}{c} -\alpha k\operatorname{col}\left\{ \frac{\partial f_{i}}{\partial x_{is}}\right\} _{i=1}^{N}\otimes e_{s}\\ \mathbf{0}_{N+Nd}\end{array}\right],
\end{array}
\]
the phase dropping out through $\cos^{2}+\sin^{2}=1$. This is the gradient term
that drives the averaged system: the averaged dynamics sees an exact gradient
rather than a residual depending on $\vartheta_{i}$.

For $i=1,\dots,N$, recalling \eqref{eq:mix_partial},
\[
\begin{array}{@{}r@{}l@{}}
\frac{\partial}{\partial t}\left(\frac{\partial f_{i}}{\partial x_{is}}\right) & {}= -\frac{\beta}{v}\xi^{-v-1}(t)\operatorname{row}\left\{ \frac{\partial^{2}f_{i}}{\partial x_{is}\partial x_{i\ell}}\right\} _{\ell=1}^{d}\tilde{x}_{fi}\\
& \quad + \operatorname{row}\left\{ \frac{\partial^{2}f_{i}}{\partial x_{is}\partial x_{i\ell}}\right\} _{\ell=1}^{d}\dot{x}^{*}(t) \\
& \quad + \frac{\partial^{2}f_{i}}{\partial x_{is}\partial\zeta}\dot{\zeta}(t).
\end{array}
\]
The Hessian $\frac{\partial^{2}f_{i}}{\partial x_{i}^{2}}$ is
bounded by the $M_{i}$-Lipschitz continuity of $\frac{\partial f_{i}}{\partial x}$
in Assumption 1 --- strong convexity supplies only the lower bound
$\frac{\partial^{2}f_{i}}{\partial x_{i}^{2}}\geq m_{i}I$, whereas an upper
bound is what is needed here. Consequently the entries
$\frac{\partial^{2}f_{i}}{\partial x_{is}\partial x_{i\ell}}$ are bounded for
all $i,s,\ell$, and
$\|\frac{\partial\left[b_{1s},b_{2s}\right]}{\partial t}\|$ is bounded for
$s=1,\dots,d$. Since
\[
\frac{\partial}{\partial\tilde{x}_{fi}}\left(\frac{\partial f_{i}}{\partial x_{is}}\right)=\xi^{-1}(t)\operatorname{col}\left\{ \frac{\partial^{2}f_{i}}{\partial x_{is}\partial x_{i\ell}}\right\} _{\ell=1}^{d},
\]
the norms $\|\frac{\partial\left[b_{1s},b_{2s}\right]}{\partial\varsigma}\|$ are
bounded for $s=1,\dots,d$ as well, which completes the proof.

For later reference we record what the bracket sum contributes to the averaged
system. Since $\frac{\partial f_{i}}{\partial x_{i}}=\sum_{s=1}^{d}\frac{\partial f_{i}}{\partial x_{is}}e_{s}$,
summing the first block over $s$ reassembles the stacked gradient,
\[
\frac{1}{2}\sum_{s=1}^{d}\left[b_{1s},b_{2s}\right]
= \left[\begin{array}{c}
-\frac{1}{2}\alpha k\frac{\partial G}{\partial\mathbf{x}}(\bar{\mathbf{x}},\zeta(t))\\
\mathbf{0}_{N+Nd}
\end{array}\right],
\]
while the $\bar{\eta}_{f}$ and $\bar{\mathbf{z}}_{f}$ blocks receive no bracket
contribution and therefore coincide with those of $b_{0}$. Substituting this
into $\dot{\bar{\varsigma}}=b_{0}(\bar{\varsigma},t)+\frac{1}{2}\sum_{s=1}^{d}\left[b_{1s},b_{2s}\right]$
yields the averaged system \eqref{averaged_system}.

\section{Proof of Lemma \ref{lem:invariant_set} }
\label{app_em:equilibrium}
Examining the third equation in system \eqref{averaged_system}
we derive
\[
\begin{array}{@{}l@{}}
\left(\mathbf{1}_{N}\otimes I_{d}\right)^{\mathrm{T}}\dot{\bar{\mathbf{z}}}_{f} \\ =\gamma\left(\mathbf{1}_{N}\otimes I_{d}\right)^{\mathrm{T}}\left(L\otimes I_{d}\right)\bar{\mathbf{x}}_{f}-\left(\mathbf{1}_{N}\otimes I_{d}\right)^{\mathrm{T}}\dot{\mathbf{z}}^{*}(t)\\
=\frac{\mathrm{d}}{\mathrm{d}t}\left(\frac{1}{2}\alpha k\left(\mathbf{1}_{N}\otimes I_{d}\right)^{\mathrm{T}}\frac{\partial G}{\partial\mathbf{x}}(\boldsymbol{1}_{N}\otimes x^{*}(t),\zeta(t))\right)\\
=\frac{\mathrm{d}}{\mathrm{d}t}\left(\frac{1}{2}\alpha k\sum_{i=1}^{N}\frac{\partial f_{i}\left(x^{*}(t),\zeta(t)\right)}{\partial x_{i}}\right)=0,
\end{array}
\]
where the first term vanishes because $\left(\mathbf{1}_{N}\otimes I_{d}\right)^{\mathrm{T}}\left(L\otimes I_{d}\right)=\left(\mathbf{1}_{N}^{\mathrm{T}}L\right)\otimes I_{d}=0$, since $\mathbf{1}_{N}^{\mathrm{T}}L=0$ for a weight-balanced graph (Assumption 3).
This implies that the quantity $\left(\mathbf{1}_{N}\otimes I_{d}\right)^{\mathrm{T}}\bar{\mathbf{z}}_{f}$
remains invariant, i.e., $\left(\mathbf{1}_{N}\otimes I_{d}\right)^{\mathrm{T}}\bar{\mathbf{z}}_{f}=\left(\mathbf{1}_{N}\otimes I_{d}\right)^{\mathrm{T}}\bar{\mathbf{z}}_{f0}$
with $\bar{z}_{f0}$ denoting the initial value of $\bar{z}_{f}$.

\section{Proof of Proposition \ref{prop:Pass}}

\label{sec:Proof-of-Proposition}
Note that in \eqref{averaged_system}, the dynamics
of $\bar{\mathbf{x}}_{f}$ and $\bar{\mathbf{z}}_{f}$ are decoupled
from $\dot{\bar{\eta}}_{f}$. Therefore, we first investigate the
stability of the following subsystem:
\begin{equation}
\begin{array}{@{}r@{}l@{}}
\dot{\bar{\mathbf{x}}}_{f}={} & \frac{\beta}{v}\xi^{-v}(t)\bar{\mathbf{x}}_{f}-\frac{1}{2}\alpha k\frac{\partial G}{\partial\mathbf{x}}(\bar{\mathbf{x}},\zeta(t))\\
 & +\left(-\mathbf{L}\bar{\mathbf{x}}_{f}-\bar{\mathbf{z}}_{f}-\mathbf{z}^{*}(t)\right)-\xi(t)\boldsymbol{1}_{N}\otimes\dot{x}^{*}(t),\\
\dot{\bar{\mathbf{z}}}_{f}={} & \gamma\mathbf{L}\bar{\mathbf{x}}_{f}-\dot{\mathbf{z}}^{*}(t).
\end{array}
\label{Lie_system}
\end{equation}
Similar to Lemma \ref{lem:invariant_set}, for any given $\epsilon\in\mathbb{R}^{d}$,
$\mathcal{S}\left(\epsilon\right)=\left\{ (\mathbf{x},\mathbf{z})|\left(\mathbf{1}_{N}\otimes I_{d}\right)^{\mathrm{T}}\mathbf{z}=\epsilon\right\} $
is a positively invariant set of system \eqref{Lie_system}. Let $Q=[r,R]\in\mathbb{R}^{N\times N}$
with $r,R$ as defined in~\eqref{eq:rR_def}, which is orthogonal since $\Vert r\Vert=1$. We use the notation
that $\mathcal{Q}=Q\otimes I_{d}$, $\mathbf{r}=r\otimes I_{d}$,
$\mathbf{R}=R\otimes I_{d}$. {Since the Laplacian matrix $L$ has a zero eigenvalue with the associated eigenvector $\boldsymbol{1}_N$, applying a change of variables (see \cite{kia2015distributed,yang2016distributed}) allows us to decouple the dynamics corresponding to the zero eigenvalue (see the first equation in \eqref{trans_system})}:
\begin{equation}
\begin{aligned} & \bar{\mathbf{x}}_{f}=\mathcal{Q}\bar{u}_{f},\quad\bar{\mathbf{z}}_{f}=\mathcal{Q}\bar{w}_{f}.\end{aligned}
\label{eq:change_variables}
\end{equation}
We partition the new variables as follows:
\begin{equation*}
     \bar{u}_{f}^{\mathrm{T}}=\left[\bar{u}_{f1}^{\mathrm{T}},\bar{u}_{f2:N}^{\mathrm{T}}\right],\quad
  \bar{w}_{f}^{\mathrm{T}}=\left[\bar{w}_{f1}^{\mathrm{T}},\bar{w}_{f{2:N}}^{\mathrm{T}}\right],
\end{equation*}
where $\bar{u}_{f1}^{\mathrm{T}},\bar{w}_{f1}^{\mathrm{T}}\in\mathbb{R}^{d}$.
In the new variables, system \eqref{Lie_system} reads as
\begin{equation}
\left\{
\begin{array}{@{}r@{}l@{}}
\dot{\bar{w}}_{f1} & {}= \mathbf{0}_{d}, \\
\dot{\bar{w}}_{f2:N} & {}= \gamma\mathbf{R^{\mathrm{T}}LR}\bar{u}_{f2:N} \\
    & \quad + \frac{1}{2}\alpha k\mathbf{R}^{\mathrm{T}}\frac{\mathrm{d}}{\mathrm{d}t}
        \left(\frac{\partial G}{\partial\mathbf{x}}(\boldsymbol{1}_{N}\otimes x^{*}(t),\zeta(t))\right), \\
\dot{\bar{u}}_{f1} & {}= \frac{\beta}{v}\xi^{-v}(t)\bar{u}_{f1}  - \frac{1}{2}\alpha k\mathbf{r}^{\mathrm{T}}g(\bar{x},\zeta(t)) \\
    & \quad - \xi(t)\mathbf{r}^{\mathrm{T}}\boldsymbol{1}_{N}\otimes\dot{x}^{*}(t), \\
\dot{\bar{u}}_{f2:N} & {}= \frac{\beta}{v}\xi^{-v}(t)\bar{u}_{f2:N}  - \frac{1}{2}\alpha k\mathbf{R}^{\mathrm{T}}g(\bar{x},\zeta(t)) \\
    & \quad - (\mathbf{R^{\mathrm{T}}LR})\bar{u}_{f2:N}  - \bar{w}_{f2:N},
\end{array}
\right.
\label{trans_system}
\end{equation}
where 
$$g(\bar{\mathbf{x}},\zeta(t))=\frac{\partial G}{\partial\mathbf{x}}(\bar{\mathbf{x}},\zeta(t))-\frac{\partial G}{\partial\mathbf{x}}(\boldsymbol{1}_{N}\otimes x^{*}(t),\zeta(t)).$$
Combining the initial condition $\overline{\mathbf{z}}_{f0}$ of Lemma \ref{lem:invariant_set} with
the first equation in \eqref{trans_system}, we observe that $\bar{w}_{f1}\equiv\mathbf{0}_{d}$
holds identically. Therefore, we only need to study the other three
equations. Denote $\psi_{1}=\left(\bar{u}_{f1}^{\mathrm{T}},\bar{u}_{f2:N}^{\mathrm{T}},\bar{w}_{f2:N}^{\mathrm{T}}\right)^{\mathrm{T}}$
and $\mathbf{P}=P\otimes I_{d}$ with
\[
P=\left(\begin{array}{ccc}
p_{11} & 0 & 0\\
\star  & p_{22}I_{N-1} & P_{2}\\
\star  & \star & P_{3}
\end{array}\right),
\]
and $\mathbf{\mathbf{P}_{2}}=P_{2}\otimes I_{d}$, $\mathbf{P}_{\mathbf{3}}=P_{3}\otimes I_{d}$.
For notational convenience, we introduce the following abbreviations:
\[
\begin{array}{@{}r@{}l@{}}
\mathcal{L} & {}= R^{\mathrm{T}}LR, \quad \bar{\mathcal{L}} = \mathcal{L}\otimes I_d, \quad \tilde{g} = \mathbf{R}^{\mathrm{T}}g(\bar{\mathbf{x}},\zeta(t)), \\
\dot{G}^{*} & {}= \frac{\mathrm{d}}{\mathrm{d}t}\left(\frac{\partial G}{\partial\mathbf{x}}(\boldsymbol{1}_{N}\otimes x^{*}(t),\zeta(t))\right).
\end{array}
\]
Consider the Lyapunov candidate
\begin{equation}\label{Lyapunov}
V=\psi_{1}^{\mathrm{T}}\mathbf{P}\psi_{1},
\end{equation}
which is positive definite because $p_{11}>0$ and the first inequality in \eqref{eq:LMI_1} makes the lower-right $2\times2$ block of $P$ positive definite. Computing $\dot{V} = 2\psi_{1}^{\mathrm{T}}\mathbf{P}\dot{\psi_{1}}$ along \eqref{trans_system} and grouping terms, we obtain
\begin{equation}
\dot{V} = \dot{V}_{\xi} + \dot{V}_{g} + \dot{V}_{L} + \dot{V}_{P_2} + \dot{V}_{\mathrm{pert}},
\label{V_decomp}
\end{equation}
where
\[
\begin{array}{@{}r@{}l@{}}
\dot{V}_{\xi} & {}= 2\frac{\beta}{v}\xi^{-v}(t)\Big(p_{11}\|\bar{u}_{f1}\|^{2} + p_{22}\|\bar{u}_{f2:N}\|^{2} \\
& \quad + \bar{w}_{f2:N}^{\mathrm{T}}\mathbf{P}_{\mathbf{2}}^{\mathrm{T}}\bar{u}_{f2:N}\Big), \\
\dot{V}_{g} & {}= -p_{11}\alpha k\bar{u}_{f1}^{\mathrm{T}}\mathbf{r^{\mathrm{T}}}g - p_{22}\alpha k\bar{u}_{f2:N}^{\mathrm{T}}\tilde{g}  - \alpha k\bar{w}_{f2:N}^{\mathrm{T}}\mathbf{P}_{\mathbf{2}}^{\mathrm{T}}\tilde{g}, \\
\dot{V}_{L} & {}= -2p_{22}\bar{u}_{f2:N}^{\mathrm{T}}\bar{\mathcal{L}}\bar{u}_{f2:N}  + 2\gamma\bar{u}_{f2:N}^{\mathrm{T}}\mathbf{P_{2}}\bar{\mathcal{L}}\bar{u}_{f2:N} \\
& \quad + 2\gamma\bar{w}_{f2:N}^{\mathrm{T}}\mathbf{P}_{\mathbf{3}}\bar{\mathcal{L}}\bar{u}_{f2:N}  - 2\bar{w}_{f2:N}^{\mathrm{T}}\mathbf{P}_{\mathbf{2}}^{\mathrm{T}}\bar{\mathcal{L}}\bar{u}_{f2:N}, \\
\dot{V}_{P_2} & {}= -2p_{22}\bar{u}_{f2:N}^{\mathrm{T}}\bar{w}_{f2:N} - 2\bar{w}_{f2:N}^{\mathrm{T}}\mathbf{P}_{\mathbf{2}}^{\mathrm{T}}\bar{w}_{f2:N}, \\
\dot{V}_{\mathrm{pert}} & {}= -2p_{11}\xi(t)\bar{u}_{f1}^{\mathrm{T}}\mathbf{r}^{\mathrm{T}}\boldsymbol{1}_{N}\otimes\dot{x}^{*}(t) \\
& \quad + \alpha k\bar{u}_{f2:N}^{\mathrm{T}}\mathbf{P_{2}}\mathbf{R}^{\mathrm{T}}\dot{G}^{*}  + \alpha k\bar{w}_{f2:N}^{\mathrm{T}}\mathbf{P_{3}}\mathbf{R}^{\mathrm{T}}\dot{G}^{*}.
\end{array}
\]
Here $g = g(\bar{\mathbf{x}},\zeta(t))$ for brevity. Since $\frac{\mathrm{d}}{\mathrm{d}t}\big(\frac{\partial f_{i}}{\partial x_{i}}(x^{*}(t),\zeta(t))\big)=\frac{\partial^{2}f_{i}}{\partial x_{i}^{2}}\dot{x}^{*}(t)+\frac{\partial^{2}f_{i}}{\partial x_{i}\partial\zeta}\dot{\zeta}(t)$ whose second derivatives are bounded by Assumption 1 on the compact range of $x^{*}(t),\zeta(t)$ that Assumption 2 grants, Assumption 2 gives a $k_{1}>0$ with $\|\dot{G}^{*}\| \leq k_{1}\xi^{c}(t)$, so the perturbation terms satisfy
\begin{equation}
\dot{V}_{\mathrm{pert}} \leq \Delta\xi^{c+1}(t)\|\mathbf{P}\|\|\psi_{1}\|
\label{varying}
\end{equation}
for some $\Delta > 0$. Let $m=\min\{m_{i}\}$ and $M=\max\{M_{i}\}$. By Lemmas \ref{lem:Lipschitz_grad}--\ref{lem:strong_convex},
\begin{equation}
\begin{array}{@{}r@{}l@{}}
\bar{\mathbf{x}}_{f}^{\mathrm{T}}g & {}\geq m\xi^{-1}(t)\|\bar{\mathbf{x}}_{f}\|^{2} \\
& {}= m\xi^{-1}(t)(\|\bar{u}_{f1}\|^{2} + \|\bar{u}_{f2:N}\|^{2}), \\
\|\tilde{g}\| & {}\leq M\xi^{-1}(t)\|\bar{\mathbf{x}}_{f}\|.
\end{array}
\label{ineq_grad_lip}
\end{equation}
Using \eqref{ineq_grad_lip}, the identity $\bar{\mathbf{x}}_{f} = \mathcal{Q}\bar{u}_{f}$ of \eqref{eq:change_variables} and $\xi(t)\geq1$, which turns the $\xi^{-1}(t)$ of \eqref{ineq_grad_lip} into $\xi^{-2}(t)$, we can bound the gradient-related terms:
\begin{equation}
\begin{array}{@{}r@{}l@{}}
\dot{V}_{g} & {}= -p_{11}\alpha k\bar{\mathbf{x}}_{f}^{\mathrm{T}}g - (p_{22}-p_{11})\alpha k\bar{u}_{f2:N}^{\mathrm{T}}\tilde{g} \\
& \quad - \alpha k\bar{w}_{f2:N}^{\mathrm{T}}\mathbf{P}_{\mathbf{2}}^{\mathrm{T}}\tilde{g} \\
& {}\leq -mp_{11}\alpha k\xi^{-2}(t)(\|\bar{u}_{f1}\|^{2} + \|\bar{u}_{f2:N}\|^{2}) \\
& \quad - (p_{22}-p_{11})\alpha k\bar{u}_{f2:N}^{\mathrm{T}}\tilde{g} - \alpha k\bar{w}_{f2:N}^{\mathrm{T}}\mathbf{P}_{\mathbf{2}}^{\mathrm{T}}\tilde{g}.
\end{array}
\label{eq:ineq_V2}
\end{equation}
Furthermore, \eqref{ineq_grad_lip} yields the following completing-the-square inequality:
\begin{equation} \label{eq:ineq_V3}
\delta M^{2}\xi^{-2}(t)(\|\bar{u}_{f1}\|^{2} + \|\bar{u}_{f2:N}\|^{2}) - \delta\|\tilde{g}\|^{2} \geq 0.
\end{equation}
Combining \eqref{V_decomp}--\eqref{eq:ineq_V3} and using $\xi^{2-v}(t) \leq 1$, which holds by \eqref{eq:vc_condition}, we obtain
\begin{equation}
\begin{array}{@{}r@{}l@{}}
\dot{V}
& {}\leq \left(\frac{2\beta}{v}p_{11} - mp_{11}\alpha k + \delta M^{2}\right)\xi^{-2}(t)\|\bar{u}_{f1}\|^{2} \\
& \quad + \left(\frac{2\beta}{v}p_{22} - mp_{11}\alpha k + \delta M^{2}\right)\xi^{-2}(t)\|\bar{u}_{f2:N}\|^{2} \\
& \quad + \frac{2\beta}{v}\xi^{-2}(t)\bar{w}_{f2:N}^{\mathrm{T}}\mathbf{P}_{\mathbf{2}}^{\mathrm{T}}\bar{u}_{f2:N} \\
& \quad - 2p_{22}\bar{u}_{f2:N}^{\mathrm{T}}\bar{\mathcal{L}}\bar{u}_{f2:N} + 2\gamma\bar{u}_{f2:N}^{\mathrm{T}}\mathbf{P_{2}}\bar{\mathcal{L}}\bar{u}_{f2:N} \\
& \quad + 2\gamma\bar{w}_{f2:N}^{\mathrm{T}}\mathbf{P}_{\mathbf{3}}\bar{\mathcal{L}}\bar{u}_{f2:N} - 2\bar{w}_{f2:N}^{\mathrm{T}}\mathbf{P}_{\mathbf{2}}^{\mathrm{T}}\bar{\mathcal{L}}\bar{u}_{f2:N} \\
& \quad - 2p_{22}\bar{u}_{f2:N}^{\mathrm{T}}\bar{w}_{f2:N} - 2\bar{w}_{f2:N}^{\mathrm{T}}\mathbf{P}_{\mathbf{2}}^{\mathrm{T}}\bar{w}_{f2:N} \\
& \quad - (p_{22}-p_{11})\alpha k\bar{u}_{f2:N}^{\mathrm{T}}\tilde{g} - \alpha k\bar{w}_{f2:N}^{\mathrm{T}}\mathbf{P}_{\mathbf{2}}^{\mathrm{T}}\tilde{g} - \delta\|\tilde{g}\|^{2} \\
& \quad + \Delta\xi^{c+1}(t)\|\mathbf{P}\|\|\psi_{1}\|.
\end{array}
\label{eq:Vdot_expanded}
\end{equation}
Define $\psi_{2} = [\bar{u}_{f2:N}^{\mathrm{T}}, \bar{w}_{f2:N}^{\mathrm{T}}, \tilde{g}^{\mathrm{T}}]^{\mathrm{T}}$. Collecting the terms in \eqref{eq:Vdot_expanded} into quadratic forms yields
\begin{equation}
\begin{array}{@{}r@{}l@{}}
\dot{V} & {}\leq \xi^{-2}(t)\psi_{1}^{\mathrm{T}}(\Phi_{1}\otimes I_{d})\psi_{1} \\
& \quad + \psi_{2}^{\mathrm{T}}(\Phi_{2}\otimes I_{d})\psi_{2} + \Delta\xi^{c+1}(t)\|\mathbf{P}\|\|\psi_{1}\|,
\end{array}
\label{V_bound}
\end{equation}
where
\begin{equation}
\Phi_{1} =
\begin{bmatrix}
\Phi_{11} & 0 & 0 \\
\star & \Phi_{12} & \frac{\beta}{v}P_{2} \\
\star & \star & -\frac{1}{2}(P_{2} + P_{2}^{\mathrm{T}})
\end{bmatrix}
\label{eq:phi_1}
\end{equation}
with $\Phi_{12} = \left(2\frac{\beta}{v}p_{22} - m p_{11} \alpha k + \delta M^{2}\right) I_{N-1}$, and $\Phi_{2}$ is given by \eqref{eq:phi_2}. The term $-2\bar{w}_{f2:N}^{\mathrm{T}}\mathbf{P}_{\mathbf{2}}^{\mathrm{T}}\bar{w}_{f2:N}$ of \eqref{eq:Vdot_expanded} is split between the $(3,3)$ block of $\Phi_{1}$ and the $(2,2)$ block of $\Phi_{2}$; this only loosens the bound, as $\Phi_{2}<0$ forces $P_{2}+P_{2}^{\mathrm{T}}>0$ and $\xi^{-2}(t)\leq1$.
We now establish convergence under $\Phi_{11}<0$ and $\Phi_{2}<0$ in \eqref{eq:LMI_1}.

Since $\Phi_{1}$ is block-diagonal between its $(1,1)$ entry and the lower-right $2\times 2$ block, we decompose
\begin{equation}
\begin{array}{@{}l@{}}
\xi^{-2}(t)\,\psi_{1}^{\mathrm{T}}(\Phi_{1}\otimes I_{d})\,\psi_{1} \\
= \xi^{-2}(t)\Phi_{11}\|\bar{u}_{f1}\|^{2} + \xi^{-2}(t)\,\chi^{\mathrm{T}}(\Phi_{1}^{\flat}\otimes I_{d})\,\chi,
\end{array}
\label{eq:phi1_split}
\end{equation}
where $\chi=(\bar{u}_{f2:N}^{\mathrm{T}},\bar{w}_{f2:N}^{\mathrm{T}})^{\mathrm{T}}$ and
\[
\Phi_{1}^{\flat}=\begin{bmatrix}
\Phi_{12} & \frac{\beta}{v}P_{2}\\
\star & -\frac{1}{2}(P_{2}+P_{2}^{\mathrm{T}})
\end{bmatrix}.
\]
Substituting \eqref{eq:phi1_split} into \eqref{V_bound} yields
\begin{equation}
\begin{array}{@{}r@{}l@{}}
\dot{V}
& {}\leq \Phi_{11}\,\xi^{-2}(t)\|\bar{u}_{f1}\|^{2} + \xi^{-2}(t)\,\chi^{\mathrm{T}}(\Phi_{1}^{\flat}\otimes I_{d})\,\chi \\
& \quad + \psi_{2}^{\mathrm{T}}(\Phi_{2}\otimes I_{d})\,\psi_{2} + \Delta\xi^{c+1}(t)\|\mathbf{P}\|\|\psi_{1}\|.
\end{array}
\label{eq:Vdot_split}
\end{equation}
Since $\psi_{2}=(\chi^{\mathrm{T}},\tilde{g}^{\mathrm{T}})^{\mathrm{T}}$ and $\Phi_{2}<0$,  the Rayleigh quotient bound gives
\[
\psi_{2}^{\mathrm{T}}(\Phi_{2}\otimes I_{d})\,\psi_{2}
\leq \lambda_{\max}(\Phi_{2})\|\psi_{2}\|^{2}
= \lambda_{\max}(\Phi_{2})\bigl(\|\chi\|^{2}+\|\tilde{g}\|^{2}\bigr).
\]
Since $\Phi_{2}<0$ is imposed in \eqref{eq:LMI_1}, every eigenvalue of $\Phi_{2}$ is negative, in particular $\lambda_{\max}(\Phi_{2})<0$, so dropping the nonnegative term $\|\tilde{g}\|^{2}$ therefore only loosens the bound:
\[
\lambda_{\max}(\Phi_{2})\bigl(\|\chi\|^{2}+\|\tilde{g}\|^{2}\bigr)
\leq \lambda_{\max}(\Phi_{2})\|\chi\|^{2}.
\]
Consequently, $\psi_{2}^{\mathrm{T}}(\Phi_{2}\otimes I_{d})\,\psi_{2} \leq \lambda_{\max}(\Phi_{2})\|\chi\|^{2}$ with $\lambda_{\max}(\Phi_{2})<0$. Combined with the $\Phi_{1}^{\flat}$ term, the total contribution from $\chi$ satisfies
\begin{equation}
\begin{array}{@{}l@{}}
\xi^{-2}(t)\,\chi^{\mathrm{T}}(\Phi_{1}^{\flat}\otimes I_{d})\,\chi + \psi_{2}^{\mathrm{T}}(\Phi_{2}\otimes I_{d})\,\psi_{2} \\
\leq \Big[\xi^{-2}(t)\lambda_{\max}(\Phi_{1}^{\flat}) + \lambda_{\max}(\Phi_{2})\Big]\|\chi\|^{2}.
\end{array}
\label{eq:eta_combined}
\end{equation}
Since $\xi^{-2}(t)\to 0$ by \eqref{eq:xi} and $\lambda_{\max}(\Phi_{2})<0$, there exists a finite time $T^{*}\geq t_{0}$ such that
\begin{equation}
\xi^{-2}(t)\lambda_{\max}(\Phi_{1}^{\flat}) + \lambda_{\max}(\Phi_{2}) \leq \frac{1}{2}\lambda_{\max}(\Phi_{2}),\quad \forall\,t\geq T^{*}.
\label{eq:Tstar_def}
\end{equation}
For $t\geq T^{*}$, substituting \eqref{eq:eta_combined}--\eqref{eq:Tstar_def} into \eqref{eq:Vdot_split} and using $\Phi_{11}<0$, we obtain
\begin{equation}
\begin{array}{@{}r@{}l@{}}
\dot{V}
& {}\leq \Phi_{11}\,\xi^{-2}(t)\|\bar{u}_{f1}\|^{2}+\frac{1}{2}\lambda_{\max}(\Phi_{2})\|\chi\|^{2} \\
& \quad + \Delta\xi^{c+1}(t)\|\mathbf{P}\|\|\psi_{1}\| \\
& {}\leq \mu\,\xi^{-2}(t)\|\psi_{1}\|^{2} + \Delta\xi^{c+1}(t)\|\mathbf{P}\|\|\psi_{1}\|,
\end{array}
\label{eq:Vdot_after_Tstar}
\end{equation}
where $\mu=\max\big\{\Phi_{11},\,\frac{1}{2}\lambda_{\max}(\Phi_{2})\big\}<0$, noting that $\frac{1}{2}\lambda_{\max}(\Phi_{2})\|\chi\|^{2}\leq\frac{1}{2}\lambda_{\max}(\Phi_{2})\,\xi^{-2}(t)\|\chi\|^{2}$ because $\xi^{-2}(t)\leq 1$ and the negative factor $\frac{1}{2}\lambda_{\max}(\Phi_{2})$ reverses this inequality, i.e.\ $\frac{1}{2}\lambda_{\max}(\Phi_{2})\left(\xi^{-2}(t)-1\right)\geq0$. By completing the square in \eqref{eq:Vdot_after_Tstar}, for $\|\psi_{1}\|\geq-\frac{2\Delta\xi^{c+3}(t)\|\mathbf{P}\|}{\mu}$, we have
\[
\dot{V}\leq \frac{\mu}{2\lambda_{\max}(P)}\,\xi^{-2}(t)\,V.
\]
By the comparison principle, for all $t\geq T^{*}$,
\[
V(t)\leq e^{\frac{\mu}{2\lambda_{\max}(P)}\int_{T^{*}}^{t}\xi^{-2}(s)\,\mathrm{d}s}\,V(T^{*}),
\]
\[
\forall\,\|\psi_{1}\|\geq-\frac{2\Delta\xi^{c+3}(t)\|\mathbf{P}\|}{\mu}.
\]
Since
\[
\int_{T^{*}}^{t}\xi^{-2}(s)\,\mathrm{d}s=\begin{cases}
\frac{1}{\beta}\ln\!\left(\frac{1+\beta(t-t_{0})}{1+\beta(T^{*}-t_{0})}\right) & v=2,\\[4pt]
\frac{v}{\beta(v-2)}\!\left[\!\left(1+\beta(t-t_{0})\right)^{\frac{v-2}{v}}\right.\\
\quad\quad\quad\left.-\left(1+\beta(T^{*}-t_{0})\right)^{\frac{v-2}{v}}\!\right] & v>2,
\end{cases}
\]
we have $\int_{T^{*}}^{t}\xi^{-2}(s)\,\mathrm{d}s\to+\infty$ as $t\to\infty$. Since $\mu<0$, it follows that $V(t)\to0$. Moreover, because $c+3<0$ by \eqref{eq:vc_condition}, the residual ball radius $-\frac{2\Delta\xi^{c+3}(t)\|\mathbf{P}\|}{\mu}\to 0$, hence $\psi_{1}\to 0$. On the finite interval $[t_{0},T^{*}]$, \eqref{eq:Vdot_split} yields $\dot{V}\leq C_{1}V + C_{2}\sqrt{V}$ for finite constants $C_{1},C_{2}>0$, so the Gronwall inequality ensures $V(t)$ remains bounded, providing the required initial value $V(T^{*})<\infty$.

Combining the above, we conclude that $\psi_{1}\to 0$, i.e., system \eqref{Lie_system} is globally asymptotically stable at the origin \emph{with respect to} $\mathcal{S}(\mathbf{0}_{d})$. Convergence is uniform in $t_{0}$, since by \eqref{eq:xi} every bound above depends on $t$ only through $t-t_{0}$. It remains to confirm the stability of the $\bar{\eta}_{f}$-subsystem in \eqref{averaged_system}. The unforced dynamics $\dot{\bar{\eta}}_{f}=\left(\frac{\beta}{v}\xi^{-v}(t)-\omega_{h}\right)\bar{\eta}_{f}$ admits the solution
\[
\bar{\eta}_{f}(t)=\xi(t)e^{-\omega_{h}(t-t_{0})}\bar{\eta}_{f}(t_{0}),
\]
which is exponentially stable. By Assumption 1, Lemmas \ref{lem:Lipschitz_grad} and \ref{lem:strong_convex} bound $h_{fi}$ above and below by $a_{i}+\frac{M_{i}}{2}\Vert\bar{x}_{i}-x^{*}\Vert^{2}$ and $a_{i}+\frac{m_{i}}{2}\Vert\bar{x}_{i}-x^{*}\Vert^{2}$, $a_{i}=\frac{\partial f_{i}}{\partial x_{i}}(x^{*},\zeta)^{\mathrm{T}}(\bar{x}_{i}-x^{*})$; as $0<m_{i}\leq M_{i}$, this bounds $\Vert h_{fi}\Vert$, giving $\xi(t)\Vert h_{fi}\Vert \leq \Vert\frac{\partial f_{i}}{\partial x_{i}}(x^{*}(t),\zeta(t))\Vert\|\bar{x}_{fi}\| + \frac{M_{i}}{2}\xi^{-1}(t)\|\bar{x}_{fi}\|^{2}$, so since $\|\frac{\partial}{\partial x_{i}}f_{i}(x^{*}(t),\zeta(t))\|$ is bounded and $\Vert\bar{\mathbf{x}}_{f}\Vert\rightarrow0$, the forcing term $\omega_{h}\xi(t)h_{f}(\bar{\mathbf{x}}_{f},\zeta(t))\rightarrow\boldsymbol{0}_{N}$. Its other forcing term vanishes too, being $O(\xi^{c+1}(t))$ with $c+1<0$ by \eqref{eq:vc_condition}. Hence the $\bar{\eta}_{f}$-system is input-to-state stable and uniformly asymptotically stable at the origin. The cascade of the two then gives Proposition \ref{prop:Pass}.


\begin{IEEEbiography}[{\includegraphics[width=1in,height=1.25in,clip,keepaspectratio]{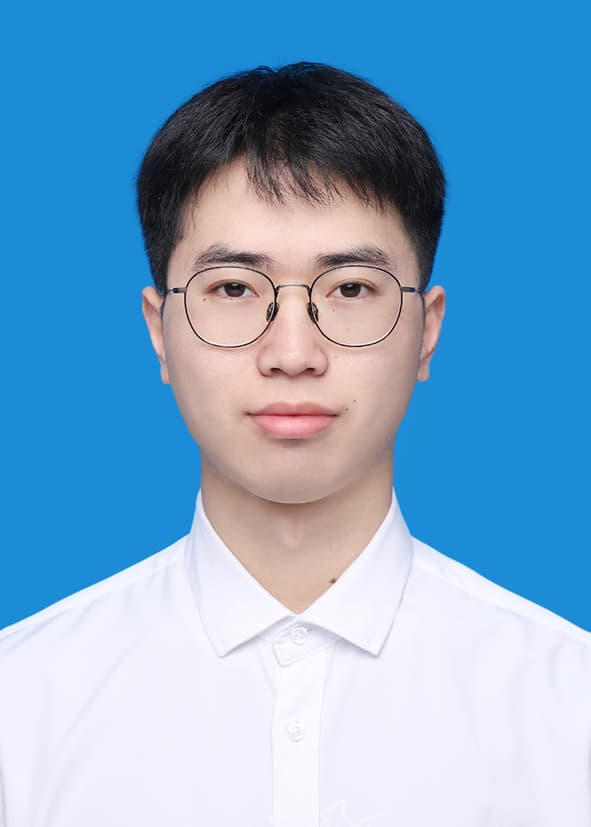}}]{Xuebin Li}
was born in 2001. He received the B.S. degree in information and computing science from the School of Mathematics, Harbin Institute of Technology, Harbin, China, in 2023. He is currently pursuing the Ph.D. degree at Harbin Institute of Technology. His research interests include optimization algorithms and reinforcement learning for robotics.
\end{IEEEbiography}


\begin{IEEEbiography}[{\includegraphics[width=1in,height=1.25in,clip,keepaspectratio]{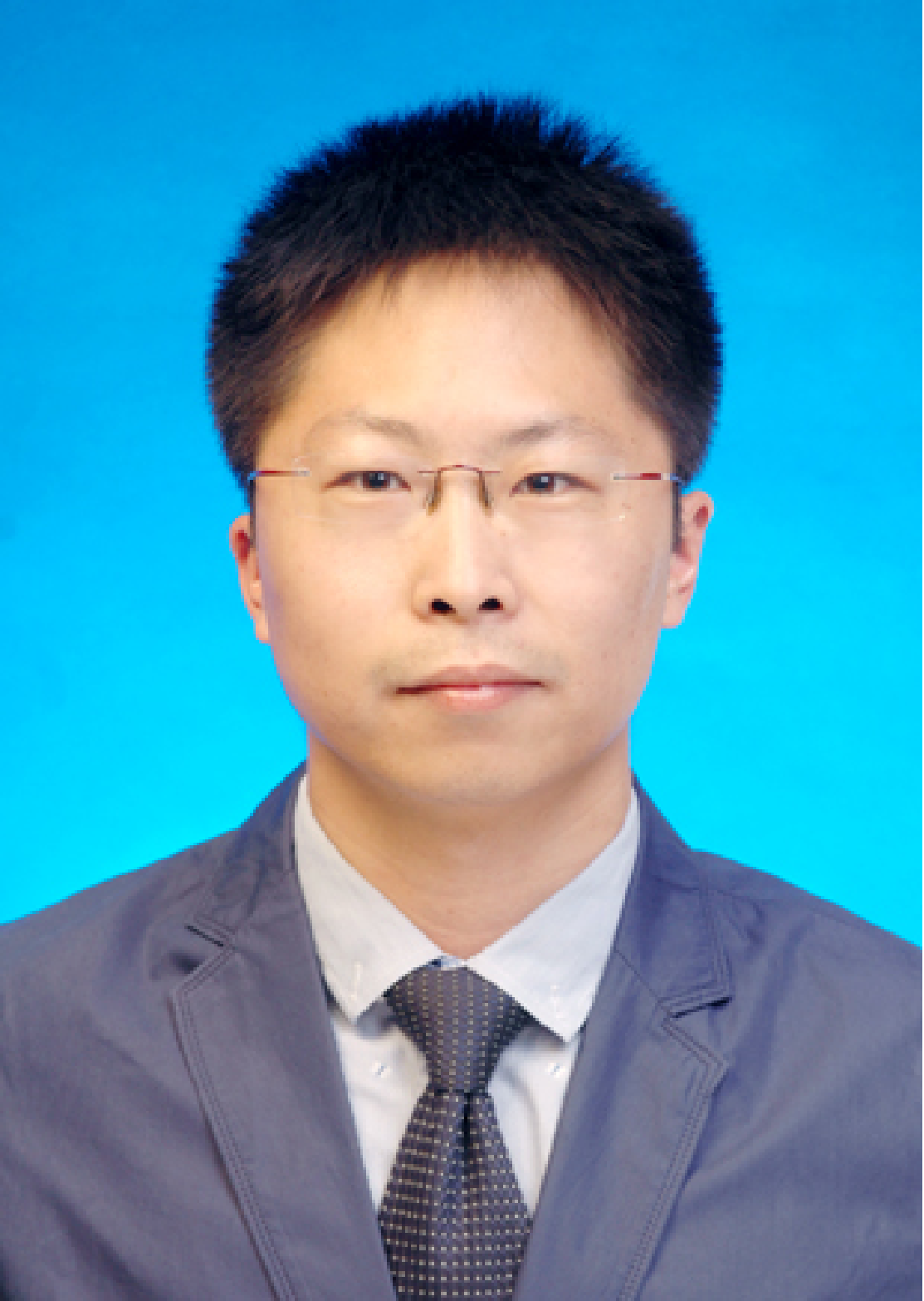}}]{Xuefei Yang}
received the B.S. degree in Mathematics, the M.S. degree in Materials Processing Engineering, the Ph.D. degree in Control Science and Engineering at Harbin Institute of Technology, Harbin, China in 2012, 2014 and 2018, respectively. From 2019 to 2021 he was a Postdoctoral Researcher at the School of Mechanical Engineering \& Automation, Harbin Institute of Technology(Shenzhen). From August 2019 to February 2020, he visited the University of Hong Kong. From 2021 to 2023 he was a Postdoctoral Researcher at the School of Electrical Engineering, Tel Aviv University, Israel. Currently he is an associate professor at Harbin Institute of Technology, Harbin. His research interests include extremum seeking control, distributed optimization and nonlinear control. He won the “National Postdoctoral Program for Innovative Talents”, the “PBC Fellowship Program for Outstanding Post-Doctoral Researchers--Israel” and the “Best Conference Paper Award--IEEE-CYBER 2023”.
\end{IEEEbiography}
\begin{IEEEbiography}[{\includegraphics[width=1in,height=1.25in,clip,keepaspectratio]{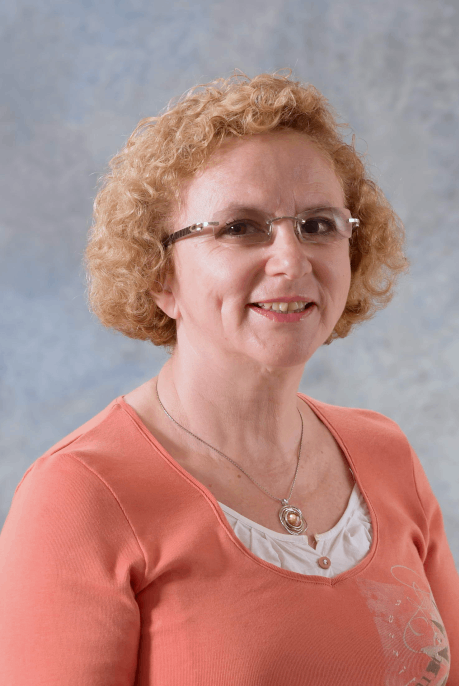}}]{Emilia Fridman} (Fellow, IEEE)
received the M.Sc. and Ph.D in mathematics in Russia. Since 1993 she has been at Tel
Aviv University, where she is currently Professor in the
Department of Electrical Engineering — Systems. She
has held numerous visiting positions in Europe, China
and Australia. Her research interests include time-delay
systems, networked control systems, distributed parameter systems, robust control and extremum seeking.
She has published more than 200 journal articles and
3 monographs. She is Editor in Chief of Systems \& Control Letters and served as Associate Editor in
Automatica, SIAM Journal on Control and Optimization
and IMA Journal of Mathematical Control and Information. She is Fellow of IEEE and IFAC. In 2014 she was ranked as a Highly Cited
Researcher by Thomson ISI. Since 2018, she has been the incumbent for Chana
and Heinrich Manderman Chair on System Control at Tel Aviv University. In 2021
she was recipient of IFAC Delay Systems Life Time Achievement Award and of
Kadar Award for outstanding research in Tel Aviv University. She was
IEEE CSS Distinguished Lecturer. In 2023 her monograph ‘Introduction to Time-Delay Systems: Analysis and Control’ (Birkhauser, 2014) was the winner of IFAC
Harold Chestnut Control Engineering Textbook Prize.
\end{IEEEbiography}
\begin{IEEEbiography}[{\includegraphics[width=1in,height=1.25in,clip,keepaspectratio]{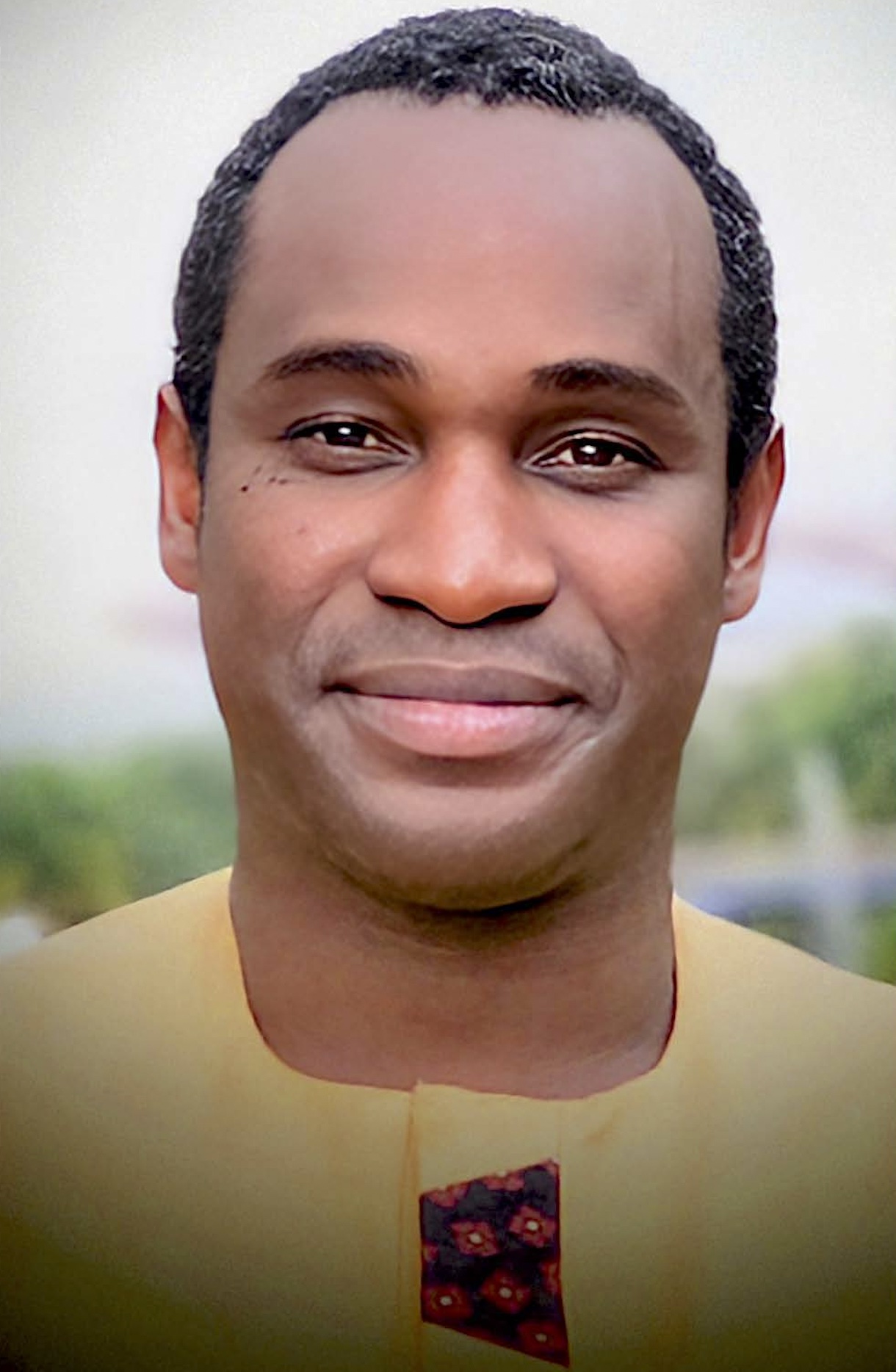}}]{Mamadou Diagne} (Senior Member, IEEE)
is with the Department of Mechanical and Aerospace Engineering, University of California San Diego, La Jolla, CA, USA, and is affiliated with the Department of Electrical and Computer Engineering. He received the Ph.D. degree in 2013 from University Claude Bernard Lyon~I, Villeurbanne, France. From 2013 to 2016, he was a Postdoctoral Scholar with the University of California San Diego and the University of Michigan. From 2017 to 2022, he was an Assistant Professor with the Department of Mechanical, Aerospace, and Nuclear Engineering, Rensselaer Polytechnic Institute, Troy, NY, USA.
His research interests include control of distributed parameter systems, adaptive control, time-delay systems, extremum seeking control, and safe nonlinear control. He served as an Associate Editor for \emph{Automatica}, \emph{Systems \& Control Letters}, and the \emph{ASME Journal of Dynamic Systems, Measurement, and Control}. He is Vice-Chair for Industry of the IFAC Technical Committee on Adaptive and Learning Systems and Vice-Chair for Education of the IFAC Technical Committee on Distributed Parameter Systems. He received the NSF CAREER Award in 2020.
\end{IEEEbiography}

\begin{IEEEbiography}[{\includegraphics[width=1in,height=1.25in,clip,keepaspectratio]{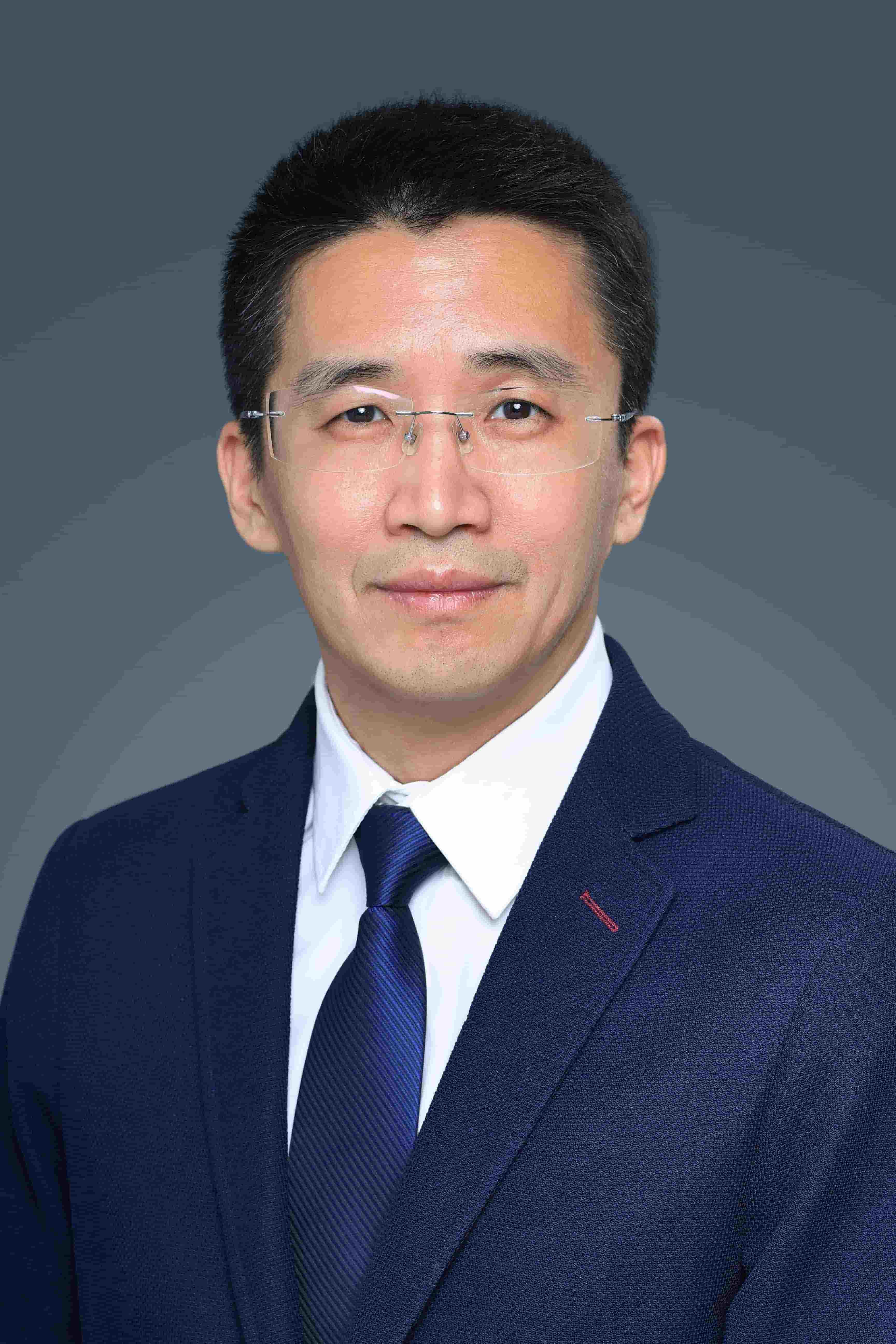}}]{Jiebao Sun}
received the Ph.D. degree from Jilin University, Changchun, China, in 2008. He then spent two years as a Postdoctoral Researcher with Harbin Institute of Technology, Harbin, China. He is currently a Professor with the School of Mathematics, Harbin Institute of Technology. His research interests include optimal control, nonlinear systems, and mathematical methods in image analysis.
\end{IEEEbiography}


\begin{thebibliography}{10}
\providecommand{\url}[1]{#1}
\csname url@samestyle\endcsname
\providecommand{\newblock}{\relax}
\providecommand{\bibinfo}[2]{#2}
\providecommand{\BIBentrySTDinterwordspacing}{\spaceskip=0pt\relax}
\providecommand{\BIBentryALTinterwordstretchfactor}{4}
\providecommand{\BIBentryALTinterwordspacing}{\spaceskip=\fontdimen2\font plus
\BIBentryALTinterwordstretchfactor\fontdimen3\font minus
  \fontdimen4\font\relax}
\providecommand{\BIBforeignlanguage}[2]{{%
\expandafter\ifx\csname l@#1\endcsname\relax
\typeout{** WARNING: IEEEtran.bst: No hyphenation pattern has been}%
\typeout{** loaded for the language `#1'. Using the pattern for}%
\typeout{** the default language instead.}%
\else
\language=\csname l@#1\endcsname
\fi
#2}}
\providecommand{\BIBdecl}{\relax}
\BIBdecl

\bibitem{ye2016distributed}
M.~Ye and G.~Hu, ``Distributed extremum seeking for constrained networked
  optimization and its application to energy consumption control in smart
  grid,'' \emph{IEEE Transactions on Control Systems Technology}, vol.~24,
  no.~6, pp. 2048--2058, 2016.

\bibitem{boyd2011distributed}
S.~Boyd, N.~Parikh, E.~Chu, B.~Peleato, J.~Eckstein \emph{et~al.},
  ``Distributed optimization and statistical learning via the alternating
  direction method of multipliers,'' \emph{Foundations and
  Trends{\textregistered} in Machine learning}, vol.~3, no.~1, pp. 1--122,
  2011.

\bibitem{de2023predefined}
P.~De~Villeros, J.~D. S{\'a}nchez-Torres, M.~Defoort, M.~Djema{\"\i}, and
  A.~Loukianov, ``Predefined-time formation control for multiagent
  systems-based on distributed optimization,'' \emph{IEEE Transactions on
  Cybernetics}, vol.~53, no.~12, pp. 7980--7988, 2023.

\bibitem{nedic2009distributed}
A.~Nedic and A.~Ozdaglar, ``Distributed subgradient methods for multi-agent
  optimization,'' \emph{IEEE Transactions on Automatic Control}, vol.~54,
  no.~1, pp. 48--61, 2009.

\bibitem{hong2017distributed}
M.~Hong, ``A distributed, asynchronous, and incremental algorithm for nonconvex
  optimization: An admm approach,'' \emph{IEEE Transactions on Control of
  Network Systems}, vol.~5, no.~3, pp. 935--945, 2017.

\bibitem{pan2024asynchronous}
Z.~Pan and M.~Cannon, ``Asynchronous admm via a data exchange server,''
  \emph{IEEE Transactions on Control of Network Systems}, vol.~11, no.~3, pp.
  1631--1643, 2024.

\bibitem{shi2015extra}
W.~Shi, Q.~Ling, G.~Wu, and W.~Yin, ``Extra: An exact first-order algorithm for
  decentralized consensus optimization,'' \emph{SIAM Journal on Optimization},
  vol.~25, no.~2, pp. 944--966, 2015.

\bibitem{brumali2026data}
R.~Brumali, G.~Carnevale, and G.~Notarstefano, ``Data-driven distributed
  optimization via aggregative tracking and deep-learning,'' \emph{IEEE
  Transactions on Control of Network Systems}, 2026.

\bibitem{wang2010control}
J.~Wang and N.~Elia, ``Control approach to distributed optimization,'' in
  \emph{2010 48th Annual Allerton Conference on Communication, Control, and
  Computing (Allerton)}.\hskip 1em plus 0.5em minus 0.4em\relax IEEE, 2010, pp.
  557--561.

\bibitem{kia2015distributed}
S.~S. Kia, J.~Cort{\'e}s, and S.~Mart{\'\i}nez, ``Distributed convex
  optimization via continuous-time coordination algorithms with discrete-time
  communication,'' \emph{Automatica}, vol.~55, pp. 254--264, 2015.

\bibitem{lu2012zero}
J.~Lu and C.~Y. Tang, ``Zero-gradient-sum algorithms for distributed convex
  optimization: The continuous-time case,'' \emph{IEEE Transactions on
  Automatic Control}, vol.~57, no.~9, pp. 2348--2354, 2012.

\bibitem{zhu2018continuous}
Y.~Zhu, W.~Yu, G.~Wen, G.~Chen, and W.~Ren, ``Continuous-time distributed
  subgradient algorithm for convex optimization with general constraints,''
  \emph{IEEE Transactions on Automatic Control}, vol.~64, no.~4, pp.
  1694--1701, 2018.

\bibitem{yang2016distributed}
S.~Yang, Q.~Liu, and J.~Wang, ``Distributed optimization based on a multiagent
  system in the presence of communication delays,'' \emph{IEEE Transactions on
  Systems, Man, and Cybernetics: Systems}, vol.~47, no.~5, pp. 717--728, 2016.

\bibitem{cherukuri2016initialization}
A.~Cherukuri and J.~Cortes, ``Initialization-free distributed coordination for
  economic dispatch under varying loads and generator commitment,''
  \emph{Automatica}, vol.~74, pp. 183--193, 2016.

\bibitem{yang2024distributed}
A.~Yang, X.~Liang, J.~Zhang, Y.~Hou, and N.~Wang, ``Distributed time-varying
  optimization with coupled constraints: Application in uav swarm
  predefined-time cooperative consensus,'' \emph{Aerospace Science and
  Technology}, vol. 147, p. 109034, 2024.

\bibitem{simonetto2016class}
A.~Simonetto, A.~Mokhtari, A.~Koppel, G.~Leus, and A.~Ribeiro, ``A class of
  prediction-correction methods for time-varying convex optimization,''
  \emph{IEEE Transactions on Signal Processing}, vol.~64, no.~17, pp.
  4576--4591, 2016.

\bibitem{fazlyab2017prediction}
M.~Fazlyab, S.~Paternain, V.~M. Preciado, and A.~Ribeiro,
  ``Prediction-correction interior-point method for time-varying convex
  optimization,'' \emph{IEEE Transactions on Automatic Control}, vol.~63,
  no.~7, pp. 1973--1986, 2017.

\bibitem{chen2025continuous}
X.~Chen, J.~I. Poveda, and N.~Li, ``Continuous-time zeroth-order dynamics with
  projection maps: Model-free feedback optimization with safety guarantees,''
  \emph{IEEE Transactions on Automatic Control}, 2025.

\bibitem{rahili2016distributed}
S.~Rahili and W.~Ren, ``Distributed continuous-time convex optimization with
  time-varying cost functions,'' \emph{IEEE Transactions on Automatic Control},
  vol.~62, no.~4, pp. 1590--1605, 2016.

\bibitem{sun2022distributed}
S.~Sun, J.~Xu, and W.~Ren, ``Distributed continuous-time algorithms for
  time-varying constrained convex optimization,'' \emph{IEEE Transactions on
  Automatic Control}, vol.~68, no.~7, pp. 3931--3946, 2022.

\bibitem{zhang2023continuous}
R.~Zhang and G.~Guo, ``Continuous distributed robust optimization of multiagent
  systems with time-varying cost,'' \emph{IEEE Transactions on Control of
  Network Systems}, vol.~11, no.~2, pp. 586--598, 2023.

\bibitem{wang2020distributed}
B.~Wang, S.~Sun, and W.~Ren, ``Distributed continuous-time algorithms for
  optimal resource allocation with time-varying quadratic cost functions,''
  \emph{IEEE Transactions on Control of Network Systems}, vol.~7, no.~4, pp.
  1974--1984, 2020.

\bibitem{wang2022distributed}
------, ``Distributed time-varying quadratic optimal resource allocation
  subject to nonidentical time-varying hessians with application to
  multiquadrotor hose transportation,'' \emph{IEEE Transactions on Systems,
  Man, and Cybernetics: Systems}, vol.~52, no.~10, pp. 6109--6119, 2022.

\bibitem{jiang2024distributed}
L.~Jiang, Z.-G. Wu, and L.~Wang, ``Distributed adaptive time-varying
  optimization with global asymptotic convergence,'' \emph{IEEE Transactions on
  Automatic Control}, 2024.

\bibitem{he2025distributed}
X.~He, Y.~Li, M.~Zhang, and T.~Huang, ``Distributed fixed-time algorithms for
  time-varying constrained optimization problems,'' \emph{IEEE Transactions on
  Artificial Intelligence}, 2025.

\bibitem{krstic2000stability}
M.~Krsti{\'c} and H.-H. Wang, ``Stability of extremum seeking feedback for
  general nonlinear dynamic systems,'' \emph{Automatica}, vol.~36, no.~4, pp.
  595--601, 2000.

\bibitem{durr2013lie}
H.-B. D{\"u}rr, M.~S. Stankovi{\'c}, C.~Ebenbauer, and K.~H. Johansson, ``Lie
  bracket approximation of extremum seeking systems,'' \emph{Automatica},
  vol.~49, no.~6, pp. 1538--1552, 2013.

\bibitem{scheinker2024100}
A.~Scheinker, ``100 years of extremum seeking: A survey,'' \emph{Automatica},
  vol. 161, p. 111481, 2024.

\bibitem{dougherty2016extremum}
S.~Dougherty and M.~Guay, ``An extremum-seeking controller for distributed
  optimization over sensor networks,'' \emph{IEEE Transactions on Automatic
  Control}, vol.~62, no.~2, pp. 928--933, 2016.

\bibitem{li2020cooperative}
Z.~Li, K.~You, and S.~Song, ``Cooperative source seeking via networked
  multi-vehicle systems,'' \emph{Automatica}, vol. 115, p. 108853, 2020.

\bibitem{mimmo2024extremum}
N.~Mimmo, G.~Carnevale, A.~Testa, and G.~Notarstefano, ``Extremum seeking
  tracking for derivative-free distributed optimization,'' \emph{IEEE
  Transactions on Control of Network Systems}, 2024.

\bibitem{nedic2017achieving}
A.~Nedic, A.~Olshevsky, and W.~Shi, ``Achieving geometric convergence for
  distributed optimization over time-varying graphs,'' \emph{SIAM Journal on
  Optimization}, vol.~27, no.~4, pp. 2597--2633, 2017.

\bibitem{yilmaz2023exponential}
C.~T. Yilmaz, M.~Diagne, and M.~Krstic, ``Exponential and prescribed-time
  extremum seeking with unbiased convergence,'' \emph{Automatica}, vol. 179, p.
  112392, 2025.

\bibitem{yilmaz2024asymptotic}
------, ``Asymptotic, exponential, and prescribed-time unbiasing in seeking of
  time-varying extrema,'' \emph{IEEE Transactions on Automatic Control}, to
  appear.

\bibitem{boyd2004convex}
S.~Boyd, ``Convex optimization,'' \emph{Cambridge UP}, 2004.

\bibitem{gharesifard2013distributed}
B.~Gharesifard and J.~Cort{\'e}s, ``Distributed continuous-time convex
  optimization on weight-balanced digraphs,'' \emph{IEEE Transactions on
  Automatic Control}, vol.~59, no.~3, pp. 781--786, 2013.

\bibitem{scheinker2014extremum}
A.~Scheinker and M.~Krsti{\'c}, ``Extremum seeking with bounded update rates,''
  \emph{Systems \& Control Letters}, vol.~63, pp. 25--31, 2014.

\bibitem{durr2013saddle}
H.-B. D{\"u}rr, C.~Zeng, and C.~Ebenbauer, ``Saddle point seeking for convex
  optimization problems,'' \emph{IFAC Proceedings Volumes}, vol.~46, no.~23,
  pp. 540--545, 2013.

\bibitem{bo2024quantization}
Y.~Bo and Y.~Wang, ``Quantization avoids saddle points in distributed
  optimization,'' \emph{Proceedings of the National Academy of Sciences}, vol.
  121, no.~18, p. e2312741121, 2024.

\bibitem{liu2018event}
S.~Liu, L.~Xie, and D.~E. Quevedo, ``Event-triggered quantized
  communication-based distributed convex optimization,'' \emph{IEEE
  Transactions on Control of Network Systems}, vol.~5, no.~1, pp. 167--178,
  2018.

\end{thebibliography}
\end{document}